\documentclass[12pt]{article}%
\usepackage{amsmath,amssymb,amsthm,amsfonts}
\usepackage{wasysym}
\usepackage{subcaption}
\usepackage{graphicx}
\usepackage[colorlinks]{hyperref}
\usepackage{color}
\usepackage{geometry}
\usepackage{appendix}
\usepackage{textcomp}

\newtheorem{thm}{Theorem}
\newtheorem{corollary}{Corollary}
\theoremstyle{definition}

\begin{document}
\title{Exotic Bifurcations Inspired by Walking Droplet Dynamics}
\author{Aminur Rahman\thanks{Corresponding Author, \url{ar276@njit.edu}}
\thanks{Department of Mathematical Sciences, New Jersey Institute of Technology}
\thanks{Current address: Department of Mathematics and Statistics, Texas Tech University} ,
Denis Blackmore\footnotemark[2]}

\date{}
\maketitle

\begin{abstract}
We identify two rather novel types of (compound)
dynamical bifurcations generated primarily by interactions of an invariant
attracting submanifold with stable and unstable manifolds of hyperbolic fixed
points. These bifurcation types - inspired by recent investigations of
mathematical models for walking droplet (pilot-wave) phenomena - are
introduced and illustrated. Some of the one-parameter bifurcation types are
analyzed in detail and extended from the plane to higher-dimensional spaces. A
few applications to walking droplet dynamics are analyzed.
\end{abstract}

\section{Introduction}

\begin{figure}[htbp]
\centering
\includegraphics[width = 0.32\textwidth]{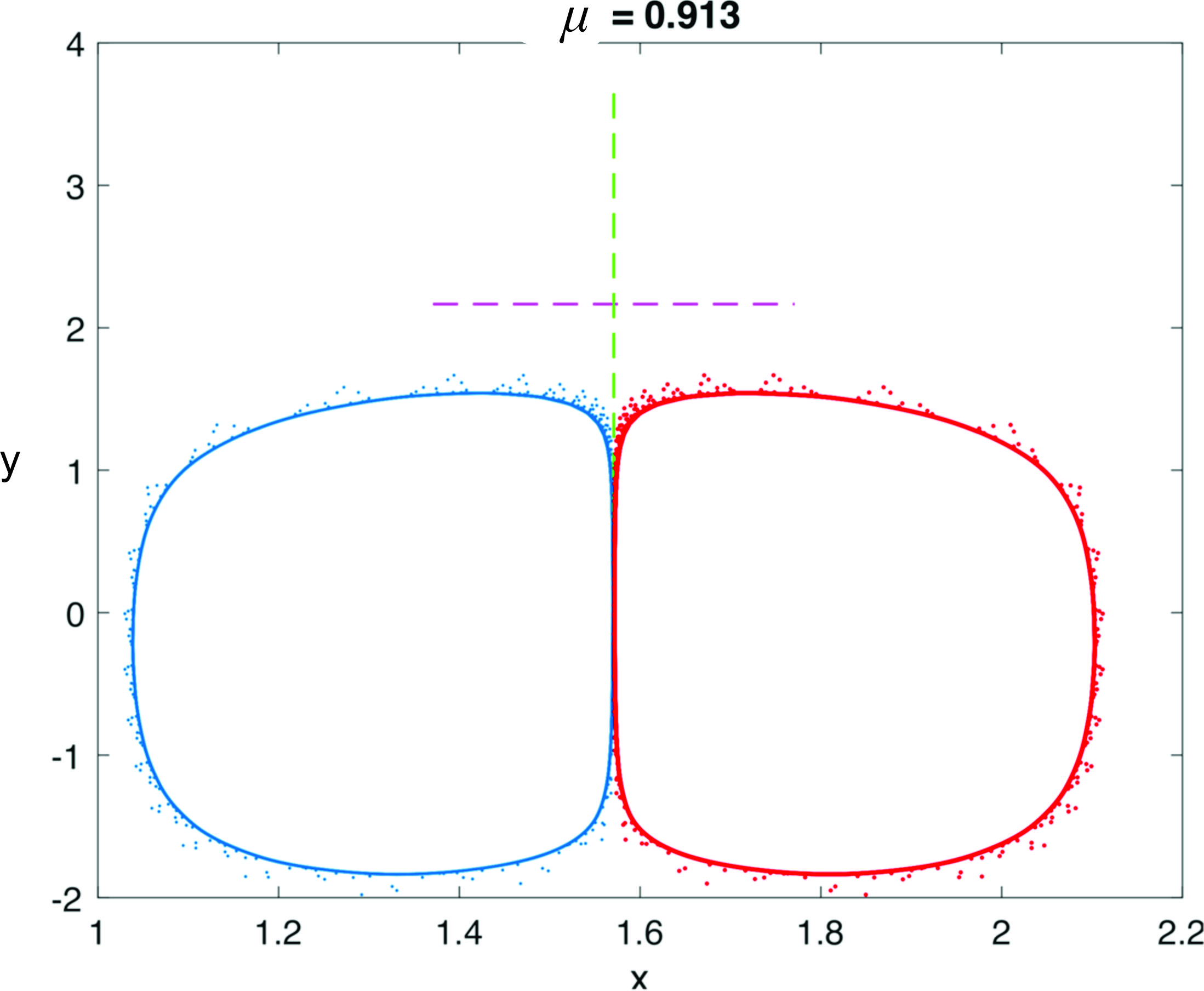}
\includegraphics[width = 0.32\textwidth]{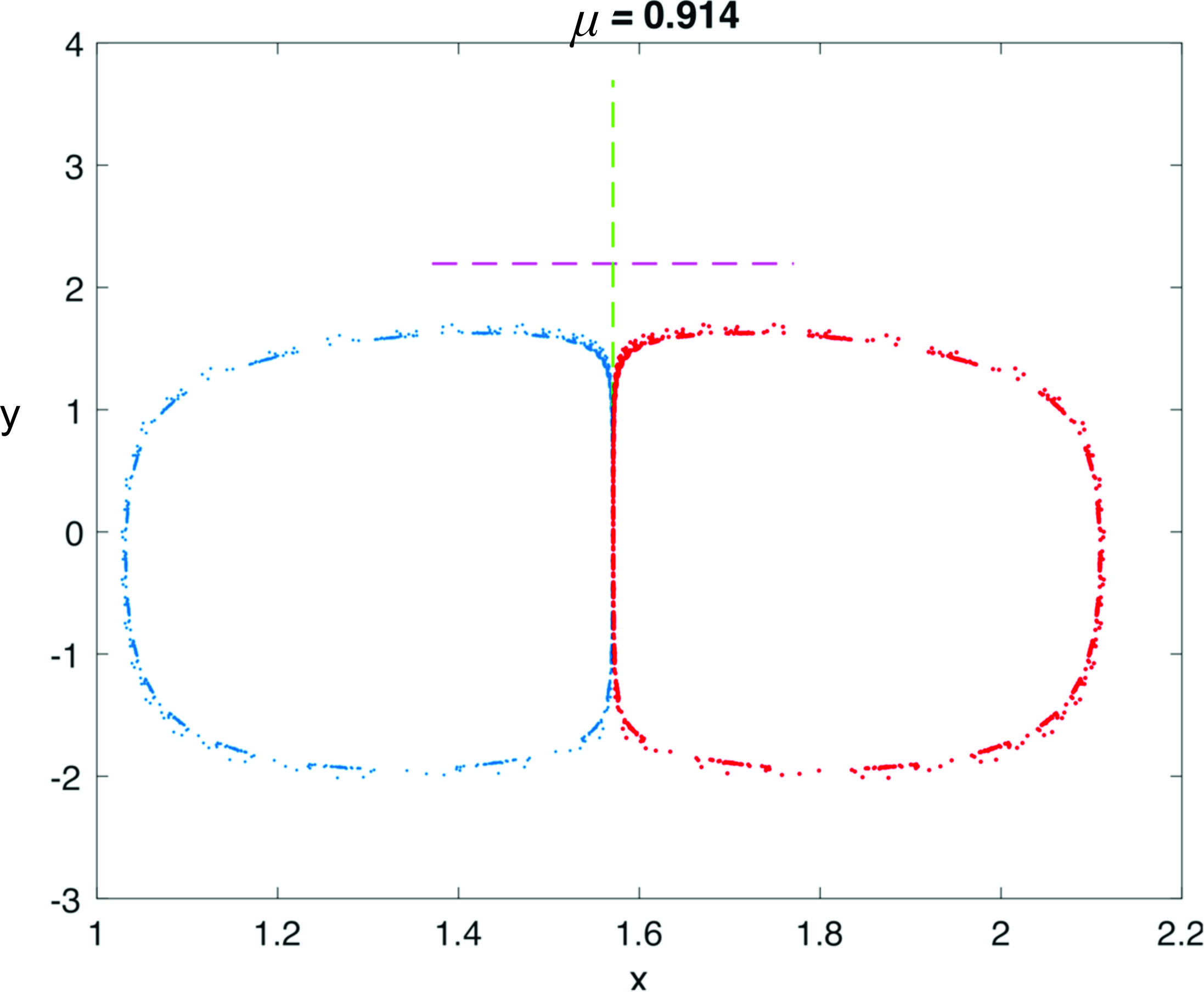}
\includegraphics[width = 0.32\textwidth]{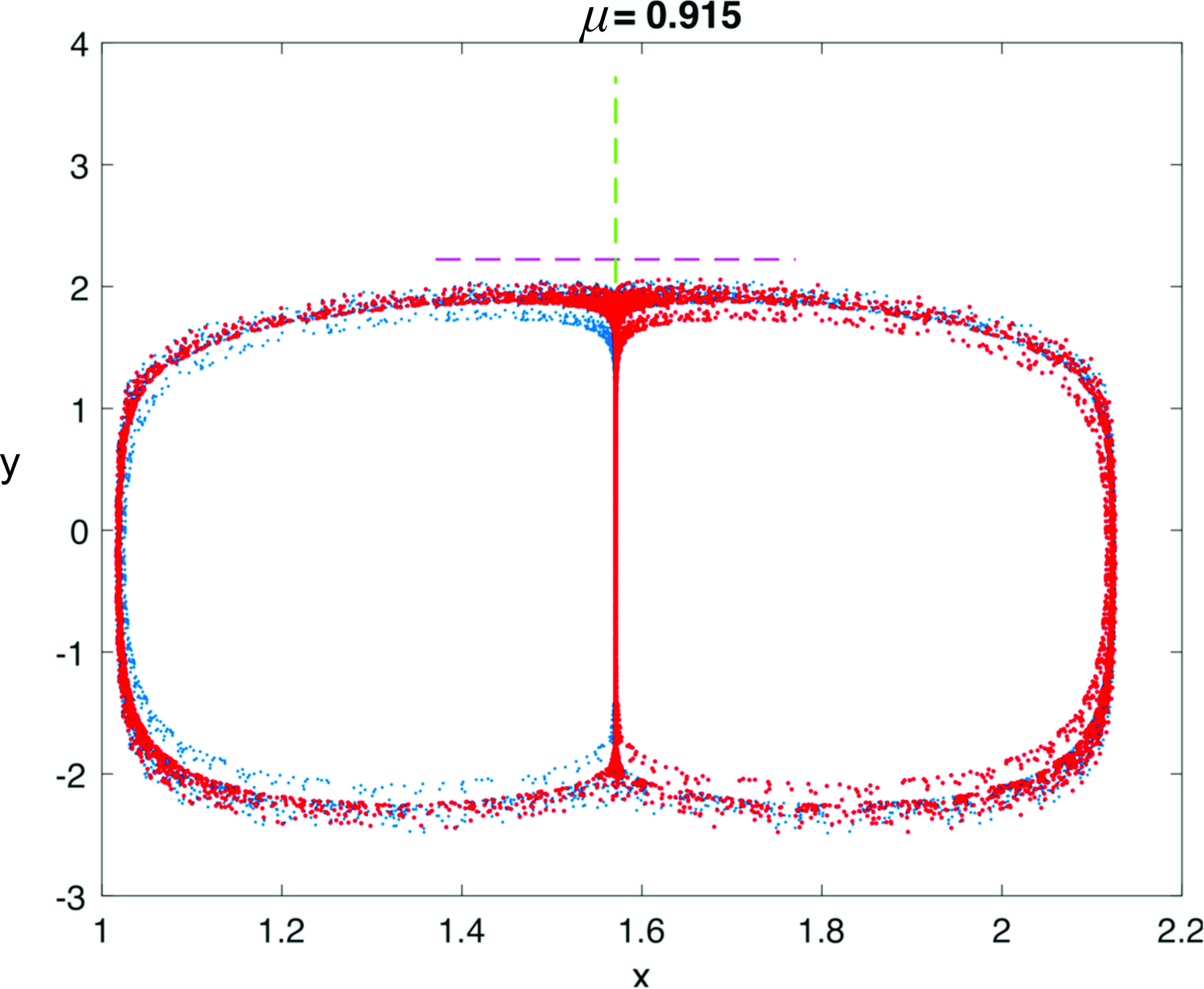}
\caption{Bifurcation evolution for fixed $C$ (increasing $\mu$).}
\label{Fig: mu-evolution}
\end{figure}

Inspired by our recent research on the dynamical properties of mathematical
models of walking droplet (pilot-wave) phenomena \cite{RB1}, we shall describe
and analyze what appear to be new types or classes of bifurcations. Owing
largely to its potential for producing macroscopic analogs of certain quantum
phenomena, walking droplet dynamics has become a very active area of research
since the seminal work of Couder \textit{et al.} \cite{CPFB}.  In this study
we focus on the dynamical systems models arising from walking droplets, interesting
examples of which can be found in Gilet \cite{Gil}, Milewski et al. \cite{MGNB}, Oza et al. \cite{OHRB},
Rahman and Blackmore \cite{RB1}, and Shirokov \cite{Shir}.  Furthermore, a detailed summary of recent
advancements in hydrodynamic pilot-waves can be found in \cite{Bush}. Simulations of the solutions of
some of the mathematical models for walkers are not only interesting for their
quantum-like effects, they exhibit exotic bifurcations that are apt to attract
the interest of dynamical systems researchers and enthusiasts. In particular,
the two-parameter planar discrete dynamical systems model of Gilet \cite{Gil}
of the form $G: \mathbb{R}^{2}\rightarrow\mathbb{R}^{2}$ defined as 
\begin{equation}
G(x,y;C;\mu):=\left(  x-C\Psi^{\prime}(x)y,\mu(y+\Psi(x))\right), \label{e1}
\end{equation}
where $0<C,\mu<1$ are parameters and $\Psi$ is an odd, $2\pi$-periodic
function given by
\[
\Psi(x):=\frac{1}{\sqrt{\pi}}\left(  \cos\beta\sin3x+\sin\beta\sin5x\right),
\]
where $\beta$ is usually chosen to be $\pi/3$ or $\pi/6$, exhibits not only 
Neimark--Sacker bifurcations, but more exotic chaotic bifurcations that have
apparently not been analyzed in detail in the literature. Simulations of the
dynamics of \eqref{e1} have shown that these exotic bifurcations are similar
in certain respects if one of the parameters $C,\mu$ is varied and the other
fixed, but also quite different in other ways. For example, in the $C$ fixed
case shown in Fig. \ref{Fig: mu-evolution} , we see a progression of similarly shaped attractors
ending in what appears to be a chaotic state. Actually, if $\mu$ is increased
further (not shown), the chaotic attractor exhibits dramatic changes, which
are apparently due to a series of dynamical crises (cf. Ott \cite{Ott}).

\begin{figure}[htbp]
\centering
\includegraphics[width = 0.49\textwidth]{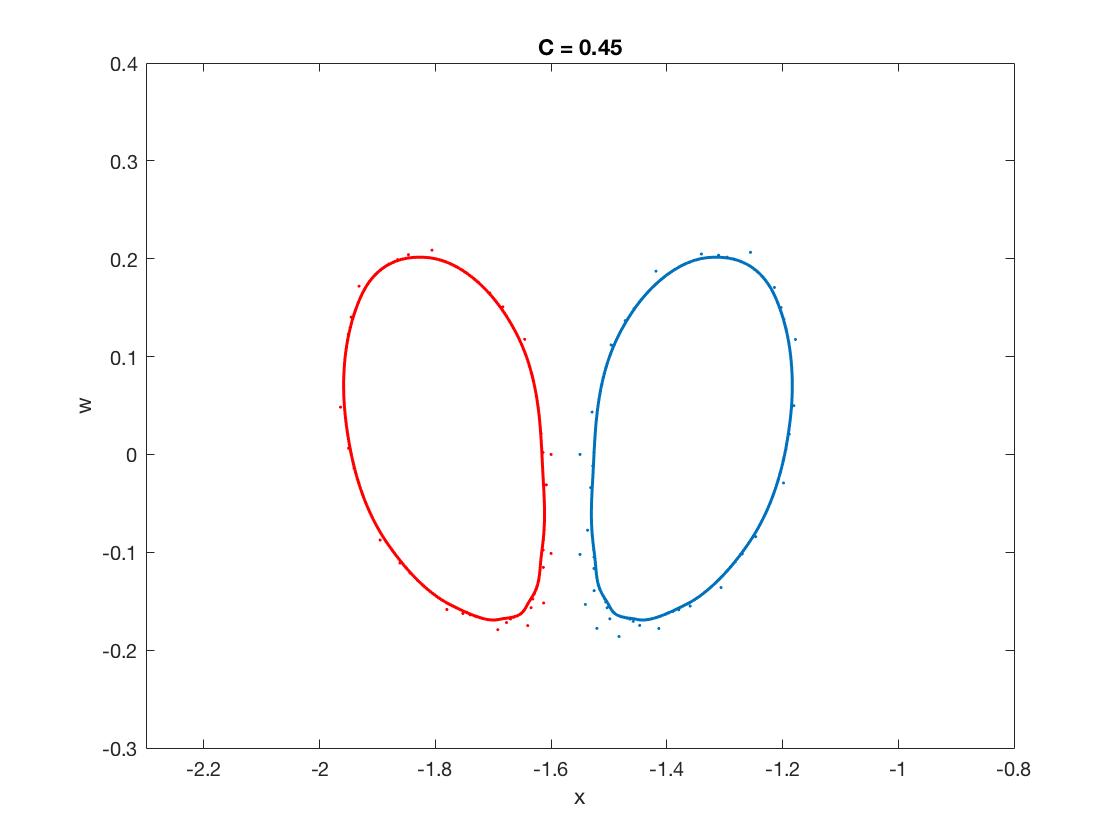}
\includegraphics[width = 0.49\textwidth]{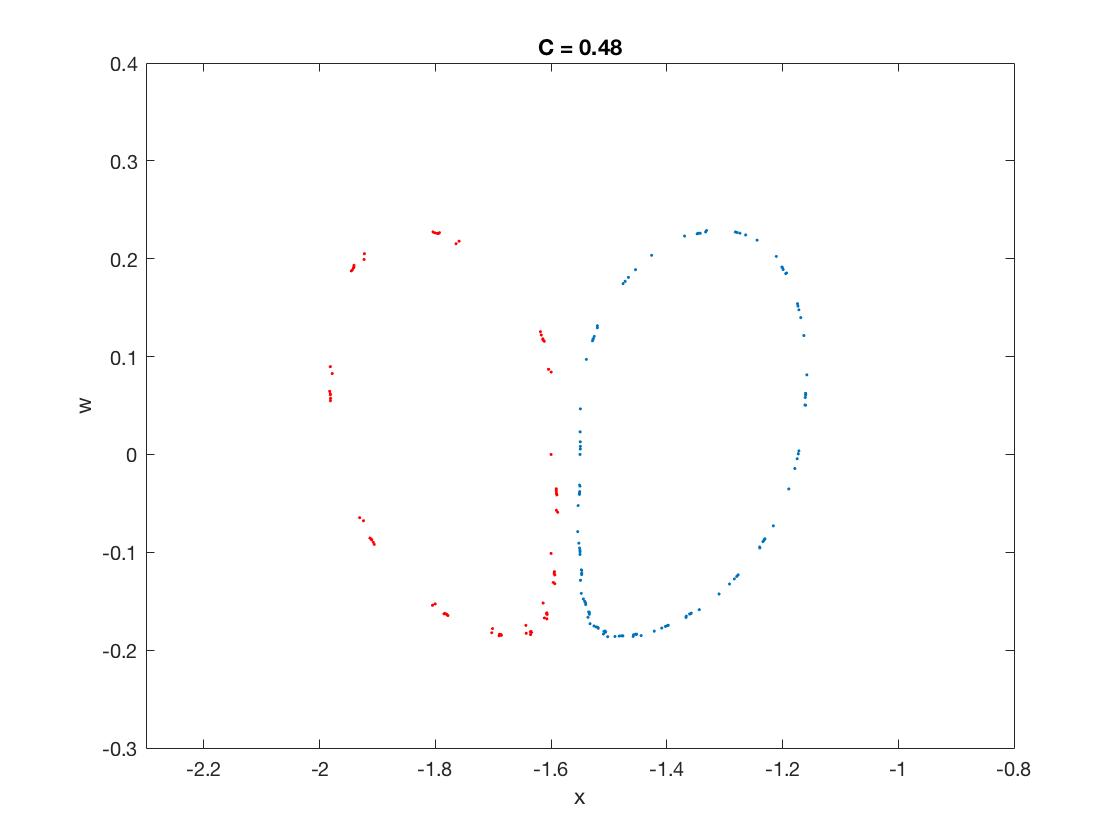}
\includegraphics[width = 0.49\textwidth]{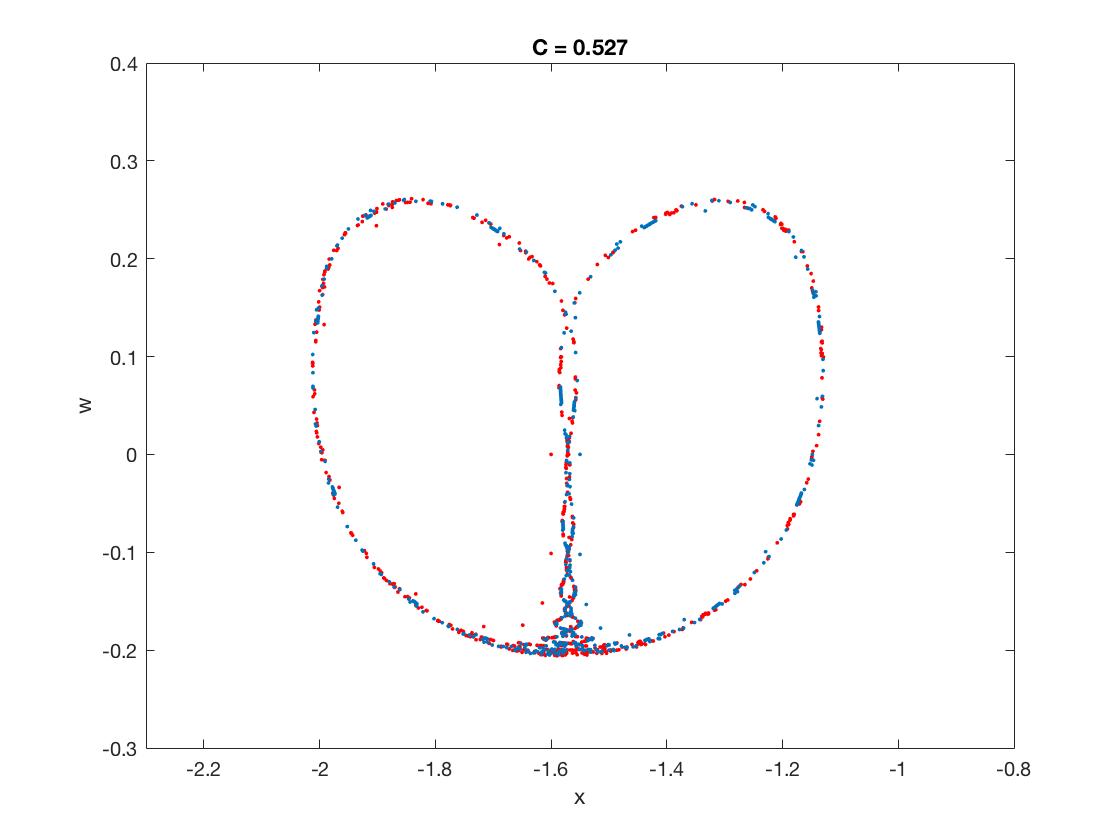}
\includegraphics[width = 0.49\textwidth]{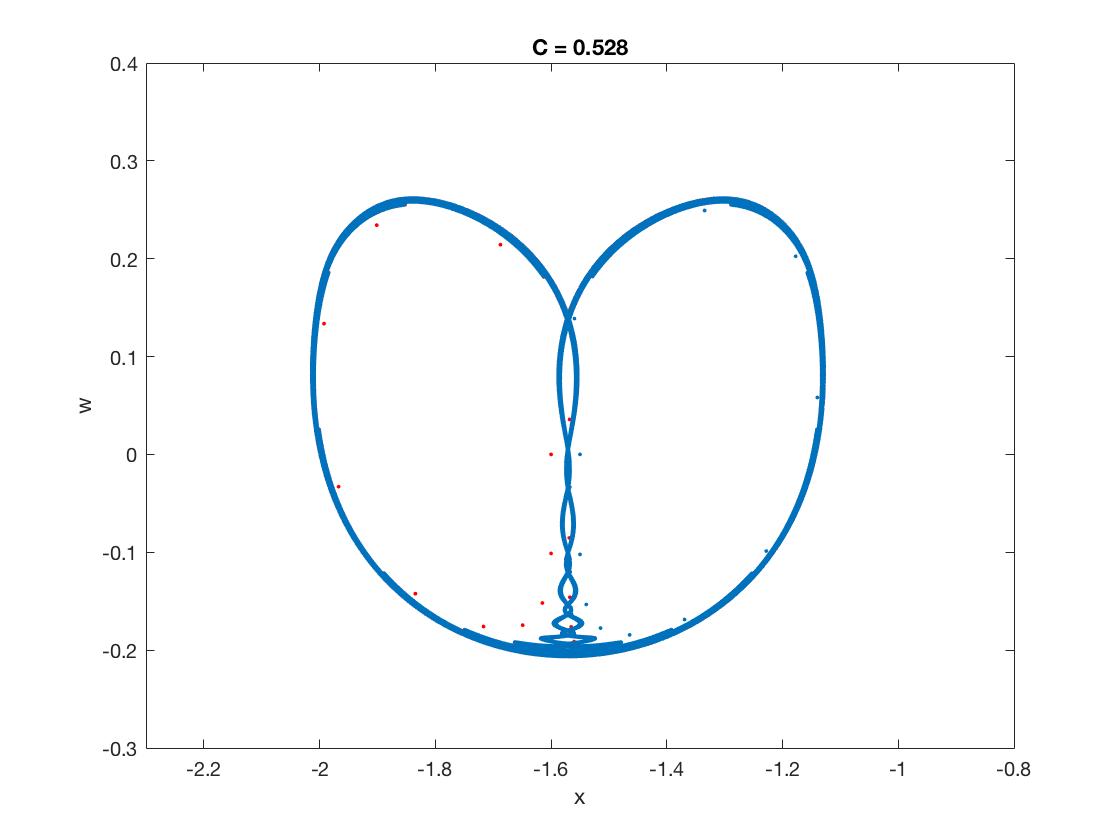}
\caption{Bifurcation evolution for fixed $\mu$ (increasing $C$).}
\label{Fig: C-evolution}
\end{figure}

On the other hand, in Fig. \ref{Fig: C-evolution} we show a sequence in which $\mu$ is fixed and
$C$ varied from $0.45$ to $0.528$ at which point we see what appears to be a
chaotic strange attractor. This chaotic strange attractor actually persists in
shape up to around $C=0.7$ (not shown). Beyond $0.7$ (also not shown) the
attractor changes in shape more or less continuously until it finally breaks
up into a chaotic splatter with just a ghost-like shadow of the prior shape -
a last stage that is also indicative of dynamical crises. In this
investigation, we concentrate on abstracting the bifurcation properties of the
Gilet map associated with fixing $\mu$ and varying $C$.


The bifurcations that we describe in the sequel are those for discrete
dynamical systems comprising the iterates of differentiable,
parameter-dependent self-maps of smooth finite-dimensional manifolds, and they
are generated by interactions of closed (positively) invariant submanifolds
with stable and unstable manifolds of saddle points as a single parameter is
varied. In the interest of simplicity and clarity, we shall confine our
attention to maps $f\in C^{1}\left(  \mathbb{R}^{m}\times\mathbb{R}%
,\mathbb{R}^{m}\right)  $, where $m$ is a natural number (in $\mathbb{N}$)
and, as usual,%
\[
C^{1}\left(  \mathbb{R}^{m}\times\mathbb{R},\mathbb{R}^{m}\right)  :=\left\{
f:\mathbb{R}^{m}\times\mathbb{R}\rightarrow\mathbb{R}^{m}:f\text{ is
continuously differentiable}\right\}  .
\]
Here, of course, $\mathbb{R}^{m}$ and $\mathbb{R}$ are, respectively,
Euclidean $m$- and $1$-space, which represent the phase and parameter spaces
of the dynamical system and points in $\mathbb{R}^{m}\times\mathbb{R}$ to be
denoted as $(\boldsymbol{x},\mathbb{\sigma})$. An example of the kind of map
we shall be investigating is that of the form of \eqref{e1} with $\mu$ a fixed
constant in $(0,1)$ and the other parameter $C$, which we denote as $\sigma$,
varying in $(0,1)$; namely,%
\begin{equation}
F(x,y;\sigma):=\left(  x-\sigma\Psi^{\prime}(x)y,\mu(y+\Psi(x))\right)  ,
\label{e2}%
\end{equation}

In what follows, we assume a knowledge of the fundamentals of modern dynamical
systems theory such as can be found in \cite{Arr,Ott,PdM,Rob,NDS}. Our main
results are detailed in the remainder of this paper, which is arranged as
follows. In Section 2 and Section 3 we describe the planar forms of the
bifurcations involving the interaction of an attracting invariant Jordan curve
and stable and unstable manifolds of saddle points. These bifurcations share
some features with the dynamical phenomena described in Aronson \emph{et al}.
\cite{ACHM} and Frouzakis \emph{et al}. \cite{FKP}, but they appear to be
essentially new. More specifically, in Section 2, we introduce and analyze
planar dynamical bifurcations generated by the interaction of an attracting
Jordan curve and the stable manifold of a saddle point. In particular, the
interaction first induces a bifurcation caused by a tangent homoclinic orbit,
which is followed by a sequence of additional tangent homoclinic orbits
interspersed with transverse homoclinic orbits as the parameter is increased.
Ultimately, however, an increase in $\sigma$ leads to a final tangent
homoclinic orbit after which there is a parameter interval on which there a
robust chaotic strange attractor, which is amenable to abstraction. As
mentioned above, there are additional types of bifurcations for larger
parameter values, which shall not be described in detail here. It should be
noted that the bifurcations considered in this paper, which shall be
designated as being of type 1, are directly related to those of the Gilet map
\eqref{e1} where a slice-like region of non-injectivity of the map plays a key
role. There is a second type - a variant of the Gilet map - where the map is a
diffeomorphism that produces analogous bifurcations, which we shall analyze in
a forthcoming paper.

Section 3 is where we describe a modification of the bifurcation in the
preceding section generated by the interaction of an attracting Jordan curve
with a pair of stable manifolds, which can induce heteroclinic cycles that
generate chaotic strange attractors followed by homoclinic bifurcations. Next,
in Section 4, we discuss some higher-dimensional analogs of the planar
bifurcations analyzed in the preceding sections. In Section 5, we describe
several applications and examples of the bifurcations, with a focus on the
phenomena and mathematical models that inspired our work on the bifurcations;
namely, walking droplet dynamics. Finally, in Section 6, we summarize some of
the conclusions reached in this research and adumbrate possible related future
work, which includes analyzing the bifurcations in the dynamics of the map
\eqref{e1} when $C$ is fixed and $\mu$ varied.

\section{Homoclinic Bifurcations of Type 1 in the Plane}

In this and the next section, we restrict our attention to a $C^{1}$ map
$f:\mathbb{R}^{2}\times\mathbb{R}\rightarrow\mathbb{R}^{2}$, with the tacit
understanding that we could have just as well considered a 1-parameter
dependent map on a simply-connected open subset of a smooth surface. The
points of $\ \mathbb{R}^{2}\times\mathbb{R}$ shall be denoted by
$(\boldsymbol{x},\sigma)=\left(  (x,y),\sigma\right)  $ and we use the
standard notation $f_{\sigma}:\mathbb{R}^{2}\rightarrow\mathbb{R}^{2}$ for the
planar map with the parameter $\sigma$ fixed at a particular value in $(0,1)$.

Let us set the stage for the\emph{ homoclinic bifurcation of type }$1$ with
more specificity. For this, we assume that $\boldsymbol{p}(\sigma):=(\hat
{x}(\sigma),\hat{y}(\sigma))$ is a saddle point of $f_{\sigma}$, with one
eigenvalue $\lambda^{u}(\sigma)>1$ and the other $0<\lambda^{s}(\sigma)<1$ for
all $\sigma\in\lbrack a,b)\subset(0,1)$ such that $\hat{x}(a),\hat{y}(a)>0$,
where $a<b$ . We may also assume that the linear unstable manifold is
horizontal when $\sigma=a$; i.e., $W_{lin}^{u}(\boldsymbol{p}(a))$ is the line
$y=\hat{y}(a)$. In addition, we assume that the saddle point, although it may
move as the parameter is varied, remains within the open rectangle
$R_{\alpha,\beta}:=(0,\hat{x}(a)+\alpha)\times(0,\hat{y}(a)+\beta)$ for some
$\alpha,\beta>0$ and all $a\leq\nu<b$, and its stable manifold $W^{s}%
(\boldsymbol{p}(\sigma))$ lies in the open vertical strip $V_{\alpha}%
:=(0,\hat{x}(a)+\alpha)\times\mathbb{R}$ for all $a\leq\sigma<1$.

A principal feature of the bifurcation is an attracting invariant closed
tubular neighborhood $T(\sigma)$ of a $C^{1}$ closed Jordan curve
$C=C(\sigma)$, which we call an \emph{invariant tube} with center $C$ for the
map. Then there is for some $\sigma$ an associated compact invariant
attracting set
\[
\mathfrak{C}\mathcal{(\sigma)}:=%
{\displaystyle\bigcap\nolimits_{n=1}^{\infty}}
f_{\sigma}^{n}\left(  T(\sigma)\right)  ,
\]
which, for an appropriately chosen center, is equal to $C(\sigma)$ for some
initial subinterval of $[a,b)$. It should be noted that such invariant circles
and tubes rather frequently arise from Neimark--Sacker bifurcations of sinks, 
especially those of the spiral variety (\emph{cf. }\cite{Neim,Sack} and also
\cite{Arr,Rob,NDS}). We also assume that for the parameter values for which
the invariant tubes exist, $T(\sigma)$ is contained in the open rectangle
$\mathcal{R}_{\alpha,\beta}:=\{(x,y):\left\vert x\right\vert <\hat
{x}(a)+\alpha,\left\vert y\right\vert <\hat{y}(a)+\beta\}$ , and even more; it
is contained in the set $Q(\sigma)$ comprising all points in $\mathcal{R}%
_{\alpha,\beta}$ to the left of $W^{s}(\boldsymbol{p}(\sigma))$ and beneath
$W_{lin}^{u}(\boldsymbol{p}(\sigma))$, respectively. Moreover, we assume that
$f_{\sigma}$ has a positive (counterclockwise) rotation number in the sense
that the iterates of any normal section of $T(\sigma)$ completely traverse the
tube in a counterclockwise manner.

Finally, there is another important feature for maps of the type represented
by Gilet's model. Namely, it is assumed that there is a $a<\sigma_{\#}<b$ such
that for all $\sigma\in\lbrack a,\sigma_{\#})$, the basin of attraction of
$\mathfrak{C}(\sigma)$ contains all points $(x,y)$ in $\mathcal{R}%
_{\alpha,\beta}$ to the left of $W^{s}(\boldsymbol{p}(\sigma))$, except for
points in a curvilinear slice $Z(\sigma)$ below the saddle point with one edge
and the image under $f_{\sigma}$ of the other edge contained in $W^{s}%
(\boldsymbol{p}(\sigma))$, having the the property that $f_{\sigma}$ flips its
interior across the stable manifold. This slice is a key agent in producing a
chaotic strange attractor as the center of the tube expands for Gilet-type maps.

\begin{figure}[htbp]
\centering
\includegraphics[width=0.9\textwidth]{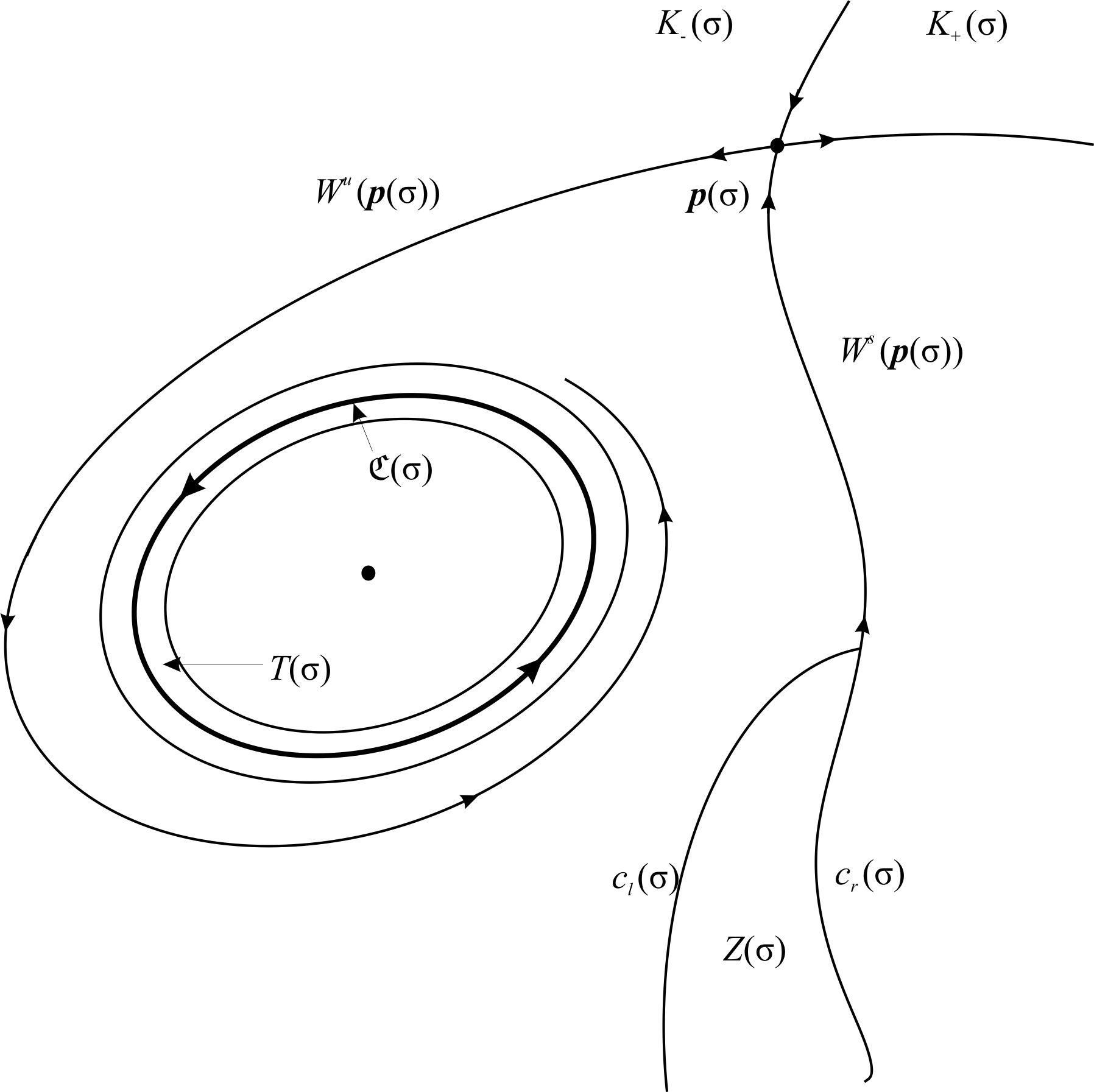}\caption{Topography of the homoclinic
type 1 bifurcation}%
\label{Top1}%
\end{figure}

Now that the contextual foundation has been established, we shall prove a
result describing the homoclinic type bifurcations that we have in mind.
Toward this end, it is useful to summarize the descriptions above in the form
of a list of detailed but reasonably succinct properties. These attributes of
the map $f:\mathbb{R}^{2}\times\mathbb{[}a,b)\rightarrow\mathbb{R}^{2}$,
illustrated in Fig. \ref{Top1}, are as follows:

\begin{itemize}
\item[(A1)] $f\in C^{1}:=C^{1}\left(  \mathbb{R}^{2}\times\mathbb{[}%
a,b),\mathbb{R}^{2}\right)  $, where $\mathbb{[}a,b)\subset(0,1)$ and
$f_{\sigma}$ has a single saddle point $\boldsymbol{p}(\sigma):=(\hat
{x}(\sigma),\hat{y}(\sigma))$, which is such that: (i) $\hat{x}(a),\hat
{y}(a)>0$; (ii) there are real constants $\kappa_{s}$, $\kappa_{u}$ such that
the eigenvalues $\lambda^{s}(\sigma)$ and $\lambda^{u}(\sigma)$ of $f_{\sigma
}^{\prime}\left(  \boldsymbol{p}(\sigma)\right)  $ satisfy $0<\lambda
^{s}(\sigma)\leq\kappa_{s}<1<\kappa_{u}\leq\lambda^{u}(\sigma)$ for all
$a\leq\sigma<b$; (iii) the eigenvector corresponding to $\lambda^{u}(a)$ is
parallel to the $x$-axis; (iv) there is a vertical strip of the form
\[
V_{\alpha}:=\{\boldsymbol{x}:=(x,y)\in\ \mathbb{R}^{2}:0<x<\hat{x}%
(a)+\alpha\}
\]
for some $\alpha>0$ such that the stable manifold $W^{s}(\boldsymbol{p}%
(\sigma))\subset S_{\alpha}$ for all $a\leq\sigma<b$ and separates the plane
into left and right components denoted as $K_{-}(\sigma)$ and $K_{+}(\sigma)$,
respectively; and (v) there is a $\beta>0$ such that
\begin{align*}
&\boldsymbol{p}(\sigma)\in R_{\alpha,\beta}:=\{(x,y)\in\
\mathbb{R}^{2}:0<x<\hat{x}(a)+\alpha,0<y<\hat{y}(a)+\beta\}\\
&\subset\mathcal{R}_{\alpha,\beta}:=\{(x,y):\left\vert x\right\vert
<\hat{x}(a)+\alpha,\left\vert y\right\vert <\hat{y}(a)+\beta\}
\end{align*}
for all $\sigma\in\lbrack a,b)$.

\item[(A2)] There exist $a<\sigma_{1}\leq\sigma_{\#}<\sigma_{t}<\sigma
_{2}<\sigma_{3}<b$ such that the following obtain: (i) there is a (positively)
$f_{\sigma}$-invariant attracting tubular neighborhood $T(\sigma)$ of a
$C^{1}$ closed Jordan curve $C(\sigma)$ for all $\sigma\in\lbrack a,\sigma
_{2})$; (ii) $f_{\sigma}$ restricted to $T(\sigma)$ has a positive
(counterclockwise) rotation number for all $\sigma\in\lbrack a,\sigma_{2})$;
(iii) the closed curve can be chosen so that \emph{centerset}
\begin{equation*}
\mathfrak{C}(\sigma):= {\displaystyle\bigcap\nolimits_{n=1}^{\infty}}
f_{\sigma}^{n}\left(  T(\sigma)\right)  =C(\sigma)
\end{equation*}
for all $\sigma\in\lbrack a,\sigma_{1})$; (iv) $T(\sigma)\subset K_{-}%
(\sigma)\cap\mathcal{R}_{\alpha,\beta}\cap H(\sigma)$, where $H(\sigma)$ is
the half-plane defined by $y<\hat{y}(\sigma)$, for all $\sigma\in\lbrack
a,\sigma_{2})$; (v) $\mathfrak{C}(\sigma)$ is a nonempty attracting set for
all $\sigma\in\lbrack a,\sigma_{\#})\cup\lbrack\sigma_{2},\sigma_{3})$ and has
an open basin of attraction, denoted as $\mathfrak{\mathring{B}}(\sigma)$,
containing $\left(  K_{-}(\sigma)\cap\mathcal{R}_{\alpha,\beta}\cap
\mathfrak{C}_{\mathrm{ext}}(\sigma)\right)  \smallsetminus Z(\sigma)$ and
$\mathfrak{C}(\sigma)\cap Z(\sigma)=\varnothing$ \ for all $\sigma\in\lbrack
a,\sigma_{\#})$, where $\mathfrak{C}_{\mathrm{ext}}(\sigma)$ is the exterior
of the centerset and $Z(\sigma)$ is a set defined as follows: There is for
every $\sigma\in\lbrack a,b)$ an orientation preserving $C^{1}$-diffeomorphism
$\Phi_{\sigma}$ of an open neighborhood $U_{\sigma}$ of a portion of
$W^{s}(\boldsymbol{p}(\sigma))$ below $\boldsymbol{p}(\sigma)$ such that
\begin{equation*}
\Phi_{\sigma}\left(  Z(\sigma)\right)  :=\left\{  (\xi,\eta):-1\leq\xi
\leq0,\eta\leq\varphi(\xi)\right\}  ,
\end{equation*}
where $\varphi:[-1,1]\rightarrow(-1,0]$ is a $C^{1}$ function such that
$\varphi(0)=0$ and $\varphi^{\prime}(\xi)>0$ when $\xi<0$. Moreover, the
right-hand bounding curve $c_{r}(\sigma):=\Phi_{\sigma}^{-1}\left(  \left\{
(0,\eta):\varphi(-1)\leq\eta\leq0\right\}  \right)  $ and the $f_{\sigma}$
image of left-hand bounding curve $c_{l}(\sigma):=\Phi_{\sigma}^{-1}\left(
\left\{  \left(  \xi,\varphi(\xi)\right)  :-1\leq\xi\leq0\right\}  \right)  $
lies in $W^{s}(\boldsymbol{p}(\sigma))$, while $f_{\sigma}$ maps the interior
$\mathring{Z}(\sigma)$ into $K_{+}(\sigma)$.
\end{itemize}

We have now sufficiently prepared the way for the statement and proof of our
first main result based on the above assumptions, which concerns the existence
of what we call a \emph{homoclinic type }$1$\emph{ bifurcation}. However, it
is convenient to first introduce the following definition: The \emph{left unstable manifold} is
\begin{equation*}
W_{l}^{u}(\boldsymbol{p}(\sigma)):=%
{\displaystyle\bigcup\nolimits_{n=1}^{\infty}}
f_{\sigma}^{n}\left(  \left\{  \boldsymbol{x}\in W^{u}(\boldsymbol{p}%
(\sigma)):\boldsymbol{x}\in\{\boldsymbol{p}(\sigma)\}\cup K_{-}(\sigma)\text{
and }d\left(  \boldsymbol{x},\boldsymbol{p}(\sigma)\right)  \leq\nu\, \forall\, \nu>0\text{ }\right\}  \right)  .
\end{equation*}

\begin{figure}[htbp]
\centering
\begin{subfigure}{0.32\textwidth}
\includegraphics[width = \textwidth]{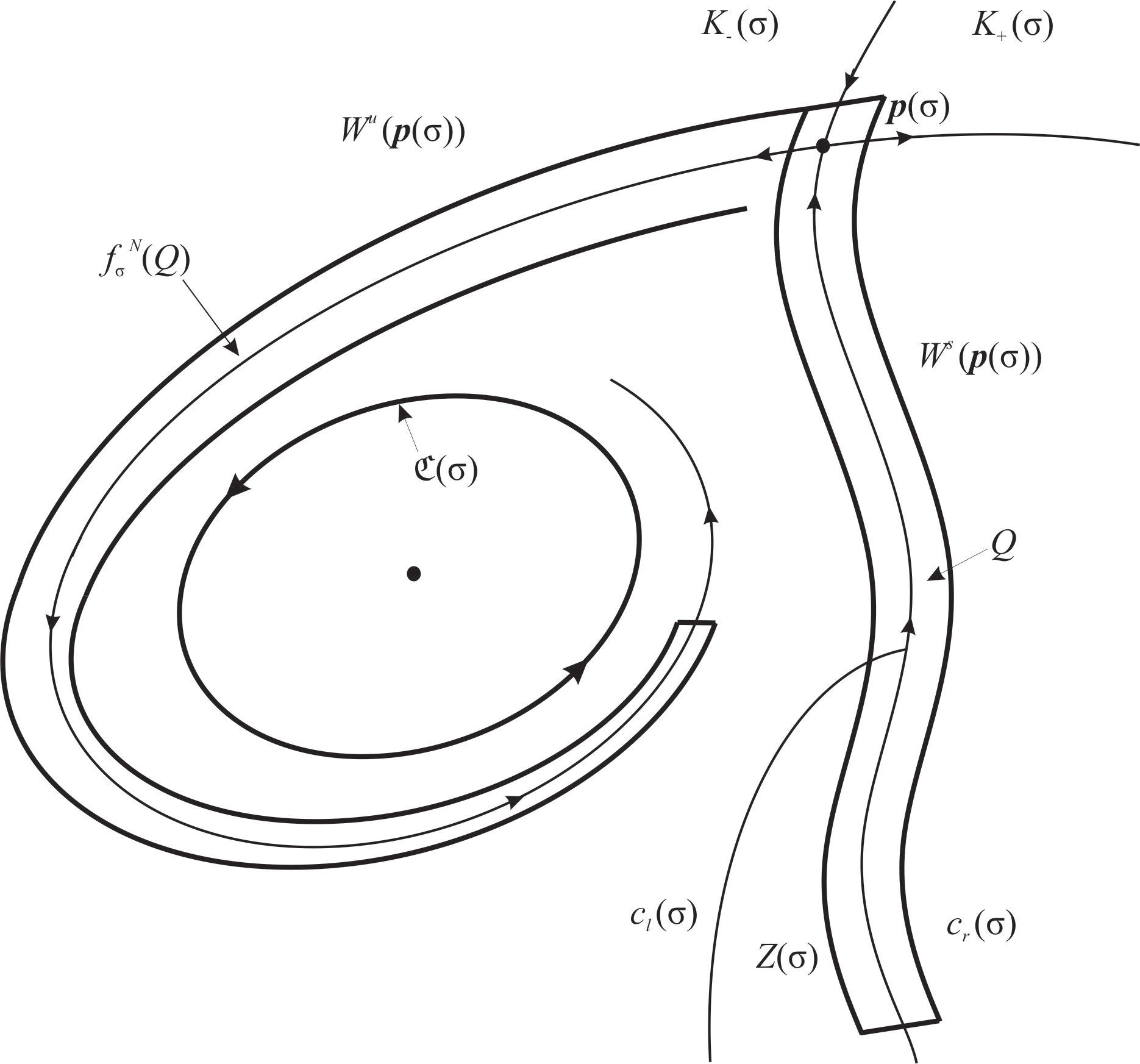}
\caption{$\sigma<\sigma_{t}$}\label{Evol-1a}
\end{subfigure}
\begin{subfigure}{0.32\textwidth}
\includegraphics[width = \textwidth]{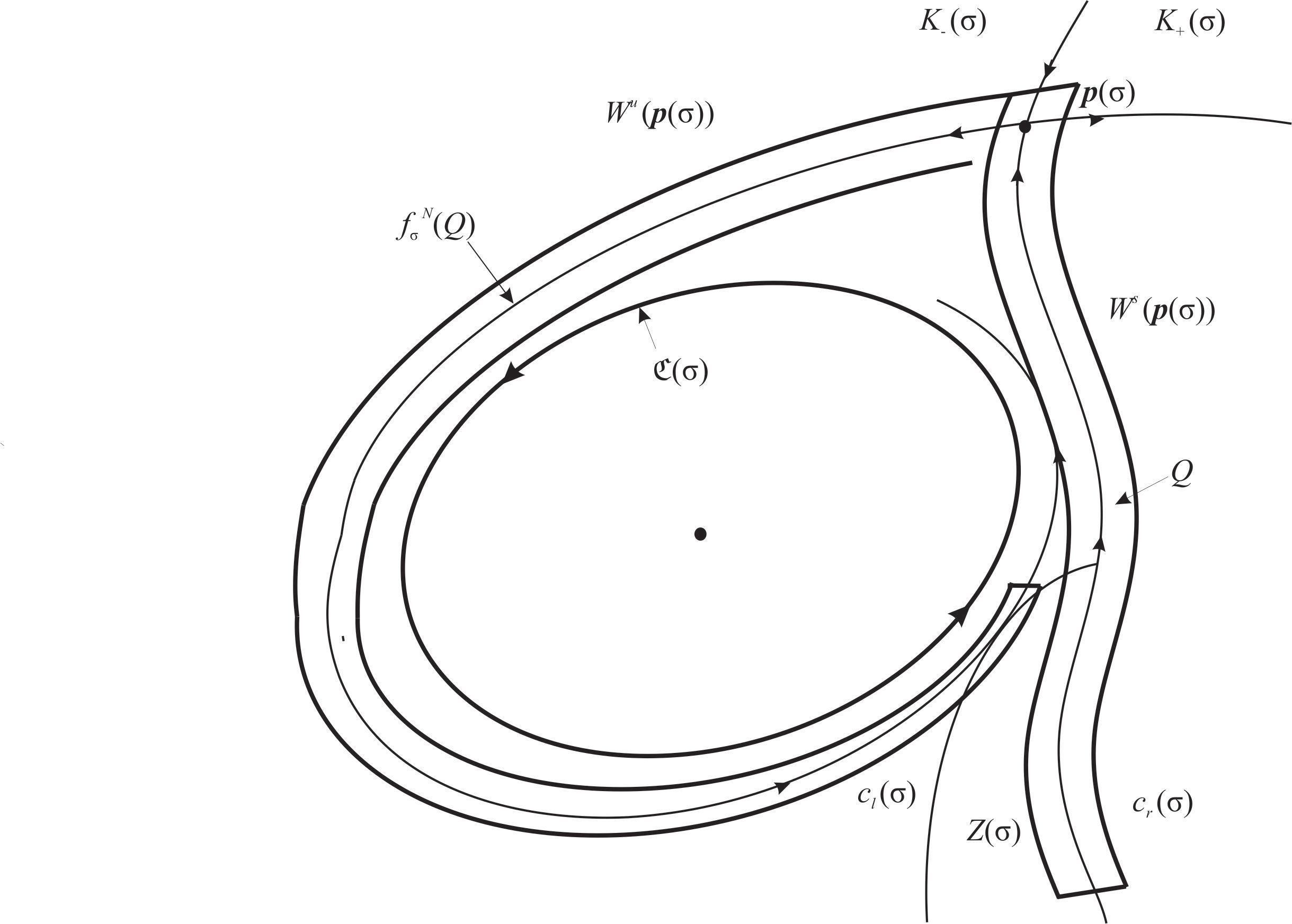}
\caption{$\sigma=\sigma_{t}$}\label{Evol-1b}
\end{subfigure}
\begin{subfigure}{0.32\textwidth}
\includegraphics[width = \textwidth]{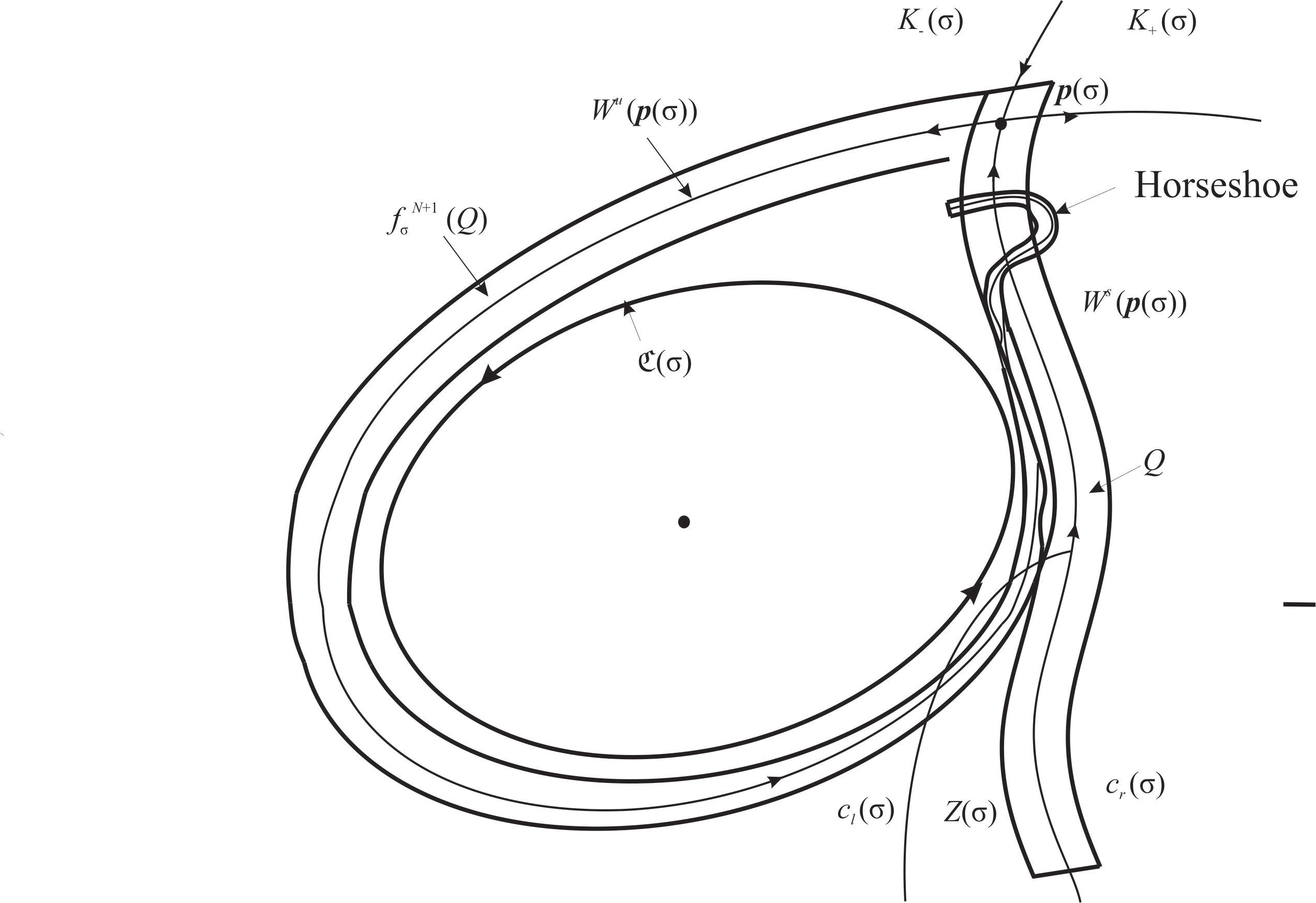}
\caption{$\sigma_{*}<\sigma<b_1$}
\label{Evol-1c}
\end{subfigure}
\caption{Geometric evolution of the homoclinic type-1 bifurcation.}
\label{Fig: Evol-1}
\end{figure}

\begin{thm}
Let $f:\mathbb{R}^{2}\times\lbrack
a,b)\rightarrow\mathbb{R}^{2}$ satisfy (A1) and (A2) and
the additional property:

\begin{itemize}
\item[(A3)] The distance $\Delta(\sigma):=\mathrm{dist}\left(
\mathfrak{C}(\sigma),Z(\sigma)\right)  :=\inf\left\{  \left\vert
\boldsymbol{x}-\boldsymbol{y}\right\vert :(\boldsymbol{x},\boldsymbol{y}%
)\in\mathfrak{C}(\sigma)\times Z(\sigma)\right\}  $ satisfies\\
$\Delta(a)>0$ and is a nonincreasing function of $\sigma$ on
$[a,\sigma_{\#})\cup\lbrack\sigma_{2},\sigma_{3})$ and $\lim
_{\sigma\uparrow\sigma_{3}}\Delta(\sigma)=0$.
\end{itemize}
Then, there are $a<\sigma_{\#}\leq\sigma_{t}\leq\sigma_{2}%
\leq\sigma_{\ast}<\sigma_{3}<b_{1}<1$, where $\sigma_{t}$ and
$\sigma_{\ast}$ are the first and last, respectively, values of
$\sigma$ where there is a tangent intersection of the stable and
unstable manifolds of $\boldsymbol{p}(\sigma)$ such that for all
$\sigma_{\ast}<\sigma<b_{1}$, $f_{\sigma}$ has a chaotic strange
attractor $\mathfrak{A}(\sigma)$, which is the closure of the left
unstable manifold: namely,
\begin{equation}
\mathfrak{A}(\sigma):=\overline{W_{l}^{u}(\boldsymbol{p}(\sigma))}\text{.}
\label{e3}%
\end{equation}

\begin{proof}
We begin by focusing on a tubular type strip
$Q=Q(\sigma;\nu)$ for $W^{s}(\boldsymbol{p}(\sigma))$ (excluding $c_{l}%
(\sigma)$) cut off just above $\boldsymbol{p}(\sigma)$ and below at the lower
edge of $R(\alpha,\beta)$ as shown in Fig. \ref{Top1}, with the understanding that the
dimensions can be taken to be as small as suits our purposes in what follows.
In particular, we take the width of the strip to be $2\nu$ and the top edge to
be parallel to and at a distance of $\nu$ above $W^{u}(\boldsymbol{p}%
(\sigma))$, where $\nu>0$. In this context, we shall find it convenient to
define $W_{\nu}^{u}(\sigma)$ to be the closed segment, with interior in
$K_{-}(\sigma)$, of the unstable manifold from $\boldsymbol{p}(\sigma)$ to the
boundary of $Q$, and to use $s_{\sigma}:=s_{\sigma}(\boldsymbol{x}%
)=s(\boldsymbol{x})$ to denote the arclength from the fixed point to any point
$\boldsymbol{x}\in W_{\nu}^{u}(\sigma)$ or any of its $f_{\sigma}$-iterates.
We also define for every nonnegative integer $n$ the \emph{endpoint}
$e_{\sigma}^{n}$, to be the boundary point of the $C^{1}$-submanifold (with
boundary) $\mathcal{W}_{\sigma}^{u}(\sigma,n):=f_{\sigma}^{n}(W_{\nu}%
^{u}(\sigma))$ for which $s_{\sigma}$ is positive and actually increases
without bound as $n\rightarrow\infty$ as long as $\sigma<\sigma_{\#}$.

The idea of the proof is illustrated rather simply in Figs. \ref{Top1} and \ref{Fig: Evol-1}, but some
technical details are necessary. It follows from the hypotheses (A1)-(A3) that
for any $\nu>0$ there is a positive integer $N=N_{\nu}(\sigma)$ and a first
$\sigma=\sigma_{t}\in\lbrack a,\sigma_{2})$ such that $\mathcal{W}_{\nu}%
^{u}(\sigma,N)$ is tangent to $c_{l}(\sigma)$, which implies that
$\mathcal{W}_{\nu}^{u}(\sigma_{t},N+1)$ is tangent to $W^{s}(\boldsymbol{p}%
(\sigma_{t}))$. Then, as $\sigma$ is increased, there must be a last value,
$\sigma_{\ast}$, with $\sigma_{t}\leq\sigma_{\ast}<\sigma_{2}$, such that
$\mathcal{W}_{\nu}^{u}(\sigma,N)$ crosses $c_{l}(\sigma)$ into the interior of
$Z(\sigma)$, implying that $Q_{l,\sigma}(\omega,N+1):=$ $f_{\sigma}%
^{N+1}(Q_{-}(\sigma;\nu))$, where $Q_{-}(\sigma;\nu):=Q(\sigma;\nu
)\cap\overline{K_{-}(\sigma)}$, which contains a portion of $\mathcal{W}_{\nu
}^{u}(\sigma,N+1)$, actually crosses over $Q(\sigma;\nu)$ into $K_{+}(\sigma)$
for all sufficiently small $\nu>0$ whenever $\sigma_{\ast}<\sigma$. As a
consequence of this situation, called an $\times$\emph{-crossing}, the desired
result follows from the attracting horseshoe theorem in \cite{JB}, which
incorporates the geometric chaos and fractal set arguments of
Birkhoff--Moser--Smale theory (\emph{cf}. \cite{Arr,Moser,PdM,Rob,Smale,NDS}).
In particular, there is an attracting horseshoe in $f_{\sigma}^{N+1}(Q)$,
which is depicted in Fig. \ref{Evol-1c}, which persists as long as $\sigma_{\ast}%
<\sigma<b_{1}$, where $\sigma=b_{1}$ is the first parameter value where the
unstable manifold has a tangent intersection with the stable manifold in
$Z(\sigma)$.

\end{proof}
\end{thm}

The proof of Theorem 1 actually yields more information than provided in the
statement; for example, it explains the \textquotedblleft blinking
effect\textquotedblright\ observed for the attractor as $\sigma:=C$ increases.
This is caused by a sequence of tangent intersections between the stable and
unstable manifolds of the saddle point $\boldsymbol{p}(\sigma)$ - each a
bifurcation value in its own right - as described by Newhouse \cite{New},
which proves the following result.

\begin{corollary}
Let $f$ be as in Theorem
$1$. Then the bifurcation value $\sigma_{\ast}$ precipitating the
creation of a robust stable chaotic strange attractor for $\sigma_{\ast
}<\sigma<b_{1}$\emph{, is proceeded by a sequence of bifurcation values
}$\sigma_{t}<\sigma_{t(1)}<\sigma_{t(2)}<\cdots<\sigma_{\ast}$
corresponding to successive tangent intersections between the stable and
unstable manifolds of the saddle point as $\sigma$ increases.
\end{corollary}

\section{Heteroclinic-Homoclinic Bifurcations of Type 1 in the Plane}

\begin{figure}[htbp]
\centering
\includegraphics[width=0.9\textwidth]{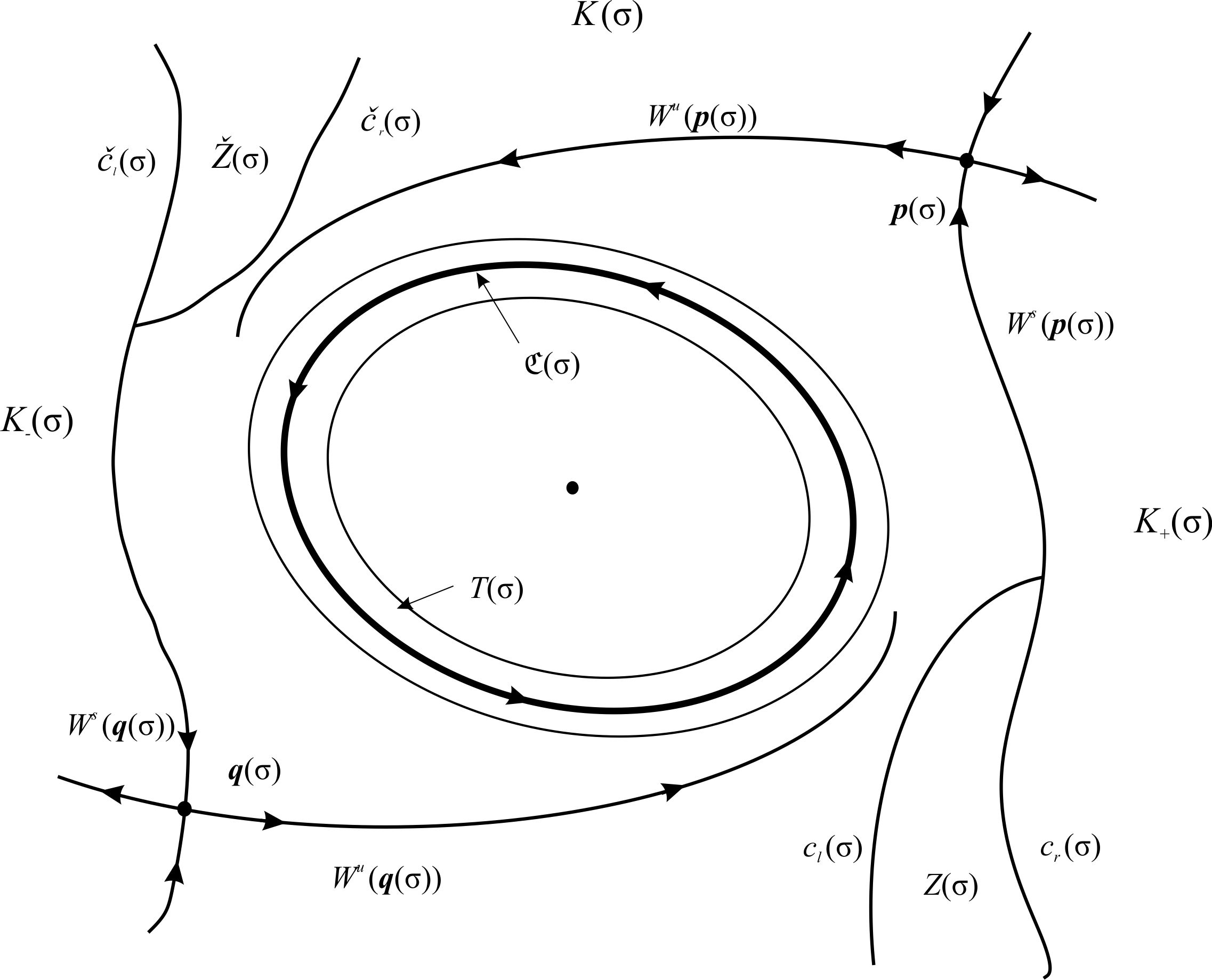}\caption{Topography of the
heteroclinic-homoclinic type 1 bifurcation}%
\label{Top2}%
\end{figure}

Next, we consider a combined heteroclinic and homoclinic variant of the
bifurcation in Section 2, illustrated in Fig. \ref{Top2}, that involves a pair of
saddle points, which we set up first with the appropriate analogs of
properties (A1) and (A2).

\begin{itemize}
\item[(B1)] $f\in C^{1}:=C^{1}\left(  \mathbb{R}^{2}\times\mathbb{[}%
a,b),\mathbb{R}^{2}\right)  $ and $f_{\sigma}$ has a pair of saddle points
$\boldsymbol{p}(\sigma):=(\hat{x}(\sigma),\hat{y}(\sigma))$ and
$\boldsymbol{q}(\sigma):=(\breve{x}(\sigma),\breve{y}(\sigma))$, which are
such that: (i) $\boldsymbol{p}(\omega)$ is in the first quadrant and
$\boldsymbol{q}(\omega)$ is in the third quadrant of $\mathbb{R}^{2}$ for all
$a\leq\sigma<b$ ; (ii) there real constants $\kappa_{s}$, $\kappa_{u}$ such
that the eigenvalues $\hat{\lambda}^{s}(\sigma)$ and $\hat{\lambda}^{u}%
(\sigma)$ of $f_{\sigma}^{\prime}\left(  \boldsymbol{p}(\sigma)\right)  $ and
$\tilde{\lambda}^{s}(\sigma)$ and $\tilde{\lambda}^{u}(\sigma)$ of $f_{\sigma
}^{\prime}\left(  \boldsymbol{q}(\omega)\right)  $ satisfy $0<\hat{\lambda
}^{s}(\sigma),\tilde{\lambda}^{s}(\sigma)\leq\kappa_{s}<1<\kappa_{u}\leq
\hat{\lambda}^{u}(\sigma),\tilde{\lambda}^{u}(\sigma)$ for all $a\leq\sigma
<b$; (iii) the eigenvectors corresponding to $\hat{\lambda}^{u}(a)$ and
$\tilde{\lambda}^{u}(a)$ are parallel to the $x$-axis; (iv) there is a
vertical strip of the form
\begin{equation*}
S_{\alpha}:=\{\boldsymbol{x}\in\ \mathbb{R}^{2}:\breve{x}(a)-\alpha<\left\vert x\right\vert <\hat{x}(a)+\alpha\}
\end{equation*}
for some $\alpha>0$ such that the stable manifolds $W^{s}(\boldsymbol{p}%
(\sigma)),W^{s}(\boldsymbol{q}(\sigma))\subset S_{\alpha}$ for all
$a\leq\sigma<b$ and separates the plane into left, middle and right components
denoted as $K_{-}(\sigma)$, $K(\sigma)$ and $K_{+}(\sigma)$, respectively; and 
(v) there is a $\beta>0$ such that
\begin{equation*}
\boldsymbol{p},\boldsymbol{q}\in R(\alpha,\beta):=\{(x,y)\in
\ \mathbb{R}^{2}:\breve{x}(a)-\alpha<\left\vert x\right\vert <\hat
{x}(a)+\alpha,\,\breve{y}(a)-\beta<\left\vert y\right\vert <\hat{y}(a)+\beta\}
\end{equation*}
for all $\sigma\in\lbrack a,b)$.

\item[(B2)] There exist $a<\sigma_{1}<\sigma_{2}<b$ such that the following
obtain: (i) there is a (positively) $f_{\sigma}$-invariant attracting tubular
neighborhood $T(\sigma)$ of a $C^{1}$ closed Jordan curve $C(\sigma)$ for all
$\sigma\in\lbrack a,\sigma_{2})$; (ii) $f_{\sigma}$ restricted to $T(\sigma)$
has a positive (counterclockwise) rotation number for all $\sigma\in\lbrack
a,\sigma_{2})$; (iii) the closed curve can be chosen so that \emph{centerset}
\begin{equation*}
\mathfrak{C}(\sigma):=
{\displaystyle\bigcap\nolimits_{n=1}^{\infty}}
f_{\omega}^{n}\left(  T(\sigma)\right)  =C(\sigma)
\end{equation*}
for all $\sigma\in\lbrack a,\sigma_{1})$; (iv) $T(\sigma)\subset K(\sigma)\cap
R(\alpha,\beta)$ for all $\sigma\in\lbrack a,\sigma_{2})$; (v) $T(\sigma)$ has
an open basin of attraction, denoted as $\mathfrak{\mathring{B}}_{T}(\sigma)$,
containing $\left(  K(\sigma)\cap R(\alpha,\beta)\cap\mathcal{C}%
_{\mathrm{ext}}(\sigma)\right)  \smallsetminus\left(  Z(\sigma)\cup\breve
{Z}(\sigma)\right)  $ for all $\sigma\in\lbrack a,\sigma_{2})$, where
$\mathcal{C}_{\mathrm{ext}}(\sigma)$ is the exterior of the centerset and
$Z(\sigma)$ and $\breve{Z}(\sigma)$ are sets defined as follows: For
$Z(\sigma)$ there is for every $\sigma\in\lbrack a,b)$ an orientation
preserving $C^{1}$-diffeomorphism $\Phi_{\sigma}$ of an open neighborhood
$W_{\sigma}$ of a portion of $W^{s}(\boldsymbol{p}(\sigma))$ below
$\boldsymbol{p}(\sigma)$ such that
\begin{equation*}
\Phi_{\sigma}\left(  Z(\sigma)\right)  :=\left\{  (\xi,\eta):-1\leq\xi
\leq0,\eta\leq\varphi(\xi)\right\}  ,
\end{equation*}
where $\varphi:[-1,1]\rightarrow(-1,0]$ is a $C^{1}$ function such that
$\varphi(0)=0$ and $\varphi^{\prime}(\xi)>0$ when $\xi<0$. Moreover, the
right-hand bounding curve $c_{r}(\sigma):=\Phi_{\sigma}^{-1}\left(  \left\{
(0,\eta):\varphi(-1)\leq\eta\leq0\right\}  \right)  $ and the $f_{\sigma}$
image of left-hand bounding curve $c_{l}(\sigma):=\Phi_{\sigma}^{-1}\left(
\left\{  \left(  \xi,\varphi(\xi)\right)  :-1\leq\xi\leq0\right\}  \right)  $
lie in $W^{s}(\boldsymbol{p}(\sigma))$, while $f_{\sigma}$ maps the interior
of $Z(\sigma)$ into $K_{+}(\sigma)$. Analogously, there is a sectorial region
$\breve{Z}(\sigma)$ (shown in Fig. \ref{Top2}) with vertex on $W^{s}(\boldsymbol{q}%
(\sigma))$ above $\boldsymbol{q}(\sigma)$ and interior above the vertex, such
that its (boundary) edges $\breve{c}_{r}(\sigma)$ and $\breve{c}_{l}(\sigma)$
lie in $K(\sigma)$ and on $W^{s}(\boldsymbol{q}(\sigma))$, respectively.
Moreover, $f_{\sigma}\left(  \breve{c}_{r}(\sigma)\right)  \subset
W^{s}(\boldsymbol{q}(\sigma))$ and $f_{\sigma}$ maps the interior of
$\breve{Z}(\sigma)$ into $K_{-}(\sigma).$
\end{itemize}

We now have assembled the basic elements needed to formulate our basic result
for heteroclinic-homoclinic bifurcations of type 1, of which there are several
variations. In the interest of keeping things as simple as possible, we impose
a symmetry requirement on the interaction of the invariant closed curve with
the slice sets $Z(\sigma)$ and $\breve{Z}(\sigma)$. For the result on the
heteroclinic-homoclinic bifurcation, it is convenient to introduce the\emph{
right unstable manifold.}
\begin{equation*}
W_{r}^{u}(\boldsymbol{q}(\sigma)):=%
{\displaystyle\bigcup\nolimits_{n=1}^{\infty}}
f_{\sigma}^{n}\left(  \left\{  \boldsymbol{x}\in W^{u}(\boldsymbol{p}%
(\sigma)):\boldsymbol{x}\in\{\boldsymbol{q}(\sigma)\}\cup K(\sigma)\text{ and
}d\left(  \boldsymbol{x},\boldsymbol{q}(\sigma)\right)  \leq\nu\,\forall\,\nu>0\text{ }\right\}  \right)  ,
\end{equation*}
which is naturally the analog of the left unstable manifold of Theorem 1 and
in the present context would be the generated by a portion the component of
the unstable manifold of $\boldsymbol{p}(\sigma)$ in $K(\sigma)$ (rather than
in what was defined as $K_{-}(\sigma)$ for Theorem 1).

\begin{thm}
Let $f:\mathbb{R}^{2}\times\lbrack
a,b)\rightarrow\mathbb{R}^{2}$satisfy (B1) and (B2) and
the additional property:

\begin{itemize}
\item[(B3)] Suppose the distance $\Delta(\sigma):=\mathrm{dist}\left(
\mathfrak{C}(\sigma),Z(\sigma)\right)  :=\inf\left\{  \left\vert
\boldsymbol{x}-\boldsymbol{y}\right\vert :(\boldsymbol{x},\boldsymbol{y}%
)\in\mathfrak{C}(\sigma)\times Z(\sigma)\right\}  =\mathrm{dist}\left(
\mathfrak{C}(\sigma),\breve{Z}(\sigma)\right)  :=\inf\left\{  \left\vert
\boldsymbol{x}-\boldsymbol{y}\right\vert :(\boldsymbol{x},\boldsymbol{y}%
)\in\mathfrak{C}(\sigma)\times\breve{Z}(\sigma)\right\}  $, satisfies
$\Delta(a)>0$ and is a nonincreasing function of $\sigma$  on
$[a,\sigma_{\#})\cup\lbrack\sigma_{2},\sigma_{3})$ and $\lim
_{\sigma\uparrow\sigma_{3}}\Delta(\sigma)=0$.
\end{itemize}
Then, there are $a<\sigma_{\#}\leq\sigma_{t}^{(1)}\leq
\sigma_{2}\leq\sigma_{\ast}^{(1)}<\sigma_{t}^{(2)}<\sigma_{\ast}^{(2)}%
<\sigma_{3}<b_{1}<1$, where $\sigma_{t}^{(1)}$and $\sigma
_{t}^{(2)}$ are the first values of $\sigma$ for which there are
tangent intersections between $W^{u}(\boldsymbol{p}(\sigma))$ and
$W^{s}(\boldsymbol{q}(\sigma))$ $($ and $W^{u}(\boldsymbol{q}%
(\sigma))$ and $W^{s}(\boldsymbol{p}(\sigma)))$ and
$W^{u}(\boldsymbol{p}(\sigma))$ and $W^{s}(\boldsymbol{p}(\sigma
))$ $($ and $W^{u}(\boldsymbol{p}(\sigma))$ and 
$W^{s}(\boldsymbol{p}(\sigma)))$, respectively. On the other hand,
$\sigma_{\ast}^{(1)}$ and $\sigma_{\ast}^{(2)}$ are the last
values of $\sigma$ where there are tangent intersections between
$W^{u}(\boldsymbol{p}(\sigma))$ and $W^{s}(\boldsymbol{q}(\sigma
))$ $($and $W^{u}(\boldsymbol{q}(\sigma))$ and %
$W^{s}(\boldsymbol{p}(\sigma)))$ and $W^{u}(\boldsymbol{p}(\sigma))$
\emph{and }$W^{s}(\boldsymbol{p}(\sigma))$ $($ and
$W^{u}(\boldsymbol{p}(\sigma))$ and $W^{s}(\boldsymbol{p}(\sigma)))$,
respectively. Then, for all $\sigma_{\ast}^{(1)}<\sigma<b_{1}$,
$f_{\sigma}$ has a chaotic strange attractor $\mathfrak{A}(\sigma
)$, which is the union of the closures of the left and right unstable
manifolds: namely,

\begin{equation}
\mathfrak{A}(\sigma):=\overline{W_{l}^{u}(\boldsymbol{p}(\sigma))}%
\cup\overline{W_{r}^{u}(\boldsymbol{q}(\sigma))}. \label{e4}%
\end{equation}

\begin{proof}
As in the proof of Theorem 1, we begin by focusing on a
tubular type strip $Q=Q(\sigma;\nu)$ for $W^{s}(\boldsymbol{p}(\sigma))$
(excluding $c_{l}(\sigma)$) cut off just above $\boldsymbol{p}(\sigma)$ and
below at the lower edge of $R(\alpha,\beta)$ (as shown in Fig. \ref{Top1}), but we also
introduce an analogous strip $\tilde{Q}=\tilde{Q}(\sigma;\nu)$ for
$W^{s}(\boldsymbol{q}(\sigma))$ with the understanding that the dimensions for
both strips can be taken to be as small as suits our purposes in what follows.
In particular, we take the width of the strips to be $2\nu$ and the top edges
to be parallel to and at a distance of $\nu$ above $W^{u}(\boldsymbol{p}%
(\sigma))$ and below $W^{u}(\boldsymbol{q}(\sigma))$, respectively, where
$\nu>0$. Again mimicking the proof of Theorem 1, but with twin objects for the
saddle points $\boldsymbol{p}$ and $\boldsymbol{q}$, it is convenient to
define $W_{\nu}^{u}(\sigma)$ and $\tilde{W}_{\nu}^{u}(\sigma)$ to be the
closed segments, with interiors in $K(\sigma)$ as shown in Fig. \ref{Top2},
respectively, of the unstable manifold from $\boldsymbol{p}(\sigma)$ to the
boundary of $Q$ and the unstable manifold from $\boldsymbol{q}(\sigma)$ to the
boundary of $\tilde{Q}$. In addition, we use $s_{\sigma}:=s_{\sigma
}(\boldsymbol{x})=s(\boldsymbol{x})$ and $\tilde{s}_{\sigma}:=\tilde
{s}_{\sigma}(\boldsymbol{x})=\tilde{s}(\boldsymbol{x})$ to denote the
arclengths from the fixed point $\boldsymbol{p}$ to any point $\boldsymbol{x}%
\in W_{\nu}^{u}(\sigma)$ and fixed point $\boldsymbol{q}$ to any point
$\boldsymbol{x}\in\tilde{W}_{\nu}^{u}(\sigma)$ or any of its $f_{\sigma}%
$-iterates. We also define for every nonnegative integer $n$ the
\emph{endpoints} $e_{\sigma}^{n}$ and $\tilde{e}_{\sigma}^{n}$ to be the
boundary points of the $C^{1}$-submanifolds (with boundary) $\mathcal{W}%
_{\sigma}^{u}(\sigma,n):=f_{\sigma}^{n}(W_{\nu}^{u}(\sigma))$ and
$\mathcal{\tilde{W}}_{\sigma}^{u}(\sigma,n):=f_{\sigma}^{n}(\tilde{W}_{\nu
}^{u}(\sigma))$ for which $s_{\sigma}$ and $\tilde{s}_{\sigma}$ are positive
and actually increase without bound as $n\rightarrow\infty$ as long as
$\sigma<\sigma_{\#}$.

As in the proof of Theorem 1, the details are best described and understood
with the aid of figures illustrating the evolution of the dynamics and
corresponding bifurcations such as in Figs. \ref{Top1} and \ref{Fig: Evol-1} for the homoclinic type
bifurcations. Although we only present the basic geometry for the case at hand
in Fig. \ref{Top2}, the analogous sequence of figures representing the evolution of
the heteroclinic-homoclinic bifurcations can easily be envisaged by comparison
with Figs. \ref{Top1} and \ref{Fig: Evol-1}, and we shall rely on this, along with an understanding of the
proof of Theorem 1, in what follows.

The hypotheses (B1)-(B3) imply that for any $\nu>0$ there are positive
integers $N^{(1)}=N_{\nu}^{(1)}(\sigma)<N^{(2)}=N_{\nu}^{(2)}(\sigma)$ such
that the following properties obtain:

\begin{itemize}
\item[(i)] There is a first $\sigma=\sigma_{t}^{(1)}\in\lbrack a,\sigma_{2})$
such that $\mathcal{W}_{\nu}^{u}(\sigma,N^{(1)})$ is tangent to $\breve{c}%
_{r}(\sigma)$ and $\mathcal{\tilde{W}}_{\nu}^{u}(\sigma,N^{(1)})$ is tangent
to $c_{l}(\sigma)$. Consequently, it follows from the definition of the slice
regions that $\mathcal{W}_{\nu}^{u}(\sigma_{t}^{(1)},N^{(1)}+1)$ is tangent to
$W^{s}(\boldsymbol{q}(\sigma_{t}^{(1)}))$ and $\mathcal{\tilde{W}}_{\nu}%
^{u}(\sigma_{t}^{(1)},N^{(1)}+1)$ is tangent to $W^{s}(\boldsymbol{p}%
(\sigma_{t}^{(1)}))$.

\item[(ii)] Moreover, there is a first $\sigma=\sigma_{t}^{(2)}\in\lbrack
a,\sigma_{2})$, with $\sigma_{t}^{(1)}<\sigma_{t}^{(2)}$, such that
$\mathcal{W}_{\nu}^{u}(\sigma,N^{(2)})$ is tangent to $c_{l}(\sigma)$ and
$\mathcal{\tilde{W}}_{\nu}^{u}(\sigma,N^{(2)})$ is tangent to $\breve{c}%
_{r}(\sigma).$ Accordingly, the characterization of the slice regions then
implies that $\mathcal{W}_{\nu}^{u}(\sigma_{t}^{(2)},N^{(2)}+1)$ is tangent to
$W^{s}(\boldsymbol{p}(\sigma_{t}^{(2)}))$ and $\mathcal{\tilde{W}}_{\nu}%
^{u}(\sigma_{t}^{(2)},N^{(2)}+1)$ is tangent to $W^{s}(\boldsymbol{q}%
(\sigma_{t}^{(2)}))$.
\end{itemize}
Therefore, we conclude from (i) that as $\sigma$ is increased, there
must be a last value, $\sigma_{\ast}^{(1)}$, with $\sigma_{t}^{(1)}\leq
\sigma_{\ast}^{(1)}<\sigma_{2}$, such that $\mathcal{W}_{\nu}^{u}%
(\sigma,N^{(1)})$ crosses $\breve{c}_{r}(\sigma)$ into the interior of
$\breve{Z}(\sigma)$, implying that $Q_{l,\sigma}(\sigma,N^{(1)}+1):=$
$f_{\sigma}^{N^{(1)}+1}(Q_{c}(\sigma;\nu))$, where $Q_{c}(\sigma
;\nu):=Q(\sigma;\nu)\cap\overline{K(\sigma)}$, which contains a portion of
$\mathcal{W}_{\nu}^{u}(\sigma,N^{(1)}+1)$, actually crosses over $\tilde
{Q}(\sigma;\nu)$ into $K_{-}(\sigma)$ for all sufficiently small $\nu>0$
whenever $\sigma_{\ast}^{(1)}<\sigma$. Furthermore, $\mathcal{\tilde{W}}_{\nu
}^{u}(\sigma,N^{(1)})$ crosses $c_{l}(\sigma)$ into the interior of
$Z(\sigma)$, implying that $\tilde{Q}_{l,\sigma}(\sigma,N^{(1)}+1):=$
$f_{\sigma}^{N^{(1)}+1}(\tilde{Q}_{c}(\sigma;\nu))$, where $\tilde{Q}%
_{c}(\sigma;\nu):=\tilde{Q}(\sigma;\nu)\cap\overline{K(\sigma)}$, which
contains a portion of $\mathcal{\tilde{W}}_{\nu}^{u}(\sigma,N^{(1)}+1)$,
actually crosses over $Q(\sigma;\nu)$ into $K_{+}(\sigma)$ for all
sufficiently small $\nu>0$ whenever $\sigma_{\ast}^{(1)}<\sigma$.
Consequently, for such values of the parameter $\sigma$ we have a transverse
heteroclinic 2-cycle connecting the saddle points $p(\sigma)$ and $q(\sigma)$,
so it follows from the work of Bertozzi \cite{Bert} (see also \cite{NDS}) that
one has a chaotic strange attractor of the form \eqref{e4} over the open
parameter interval $(\sigma_{\ast}^{(1)},\sigma_{t}^{(2)}).$

When the parameter value increases to $\sigma_{t}^{(2)}$, it follows from
(ii) that we get a small Newhouse type bifurcation \cite{New}, which amounts
to a "blinking effect" followed by possibly more tangent intersections of
$W^{u}(\boldsymbol{p}(\sigma))$ and $W^{s}(\boldsymbol{p}(\sigma))$ and
$W^{u}(\boldsymbol{q}(\sigma))$ and $W^{s}(\boldsymbol{q}(\sigma))$. Moreover,
there is a last such tangent intersection of these unstable and stable
manifolds, beyond which we have the the geometric chaos and fractal set
arguments of Birkhoff--Moser--Smale theory that imply the existence of the
chaotic strange attractor defined by \eqref{e4}, amounting to a two-fold
version of the result for a single saddle point in Theorem 1. In particular,
there is a two-fold manifestation of the attracting horseshoe in $f_{\sigma
}^{N+1}(Q)$ depicted in Fig. \ref{Evol-1c}, with a copy at each of the saddle points. In
a manner completely analogous to the homoclinic bifurcation, this (double)
horseshoe bifiurcation persists as long as $\sigma_{\ast}<\sigma<b_{1}$, where
$\sigma=b_{1}$ is the first parameter value where the unstable manifolds
$W^{u}(\boldsymbol{p}(\sigma))$ and $W^{u}(\boldsymbol{q}(\sigma))$ have
tangent intersections with the stable manifolds $W^{s}(\boldsymbol{p}%
(\sigma))$ and $W^{s}(\boldsymbol{q}(\sigma))$ in $Z(\sigma)$ and $\breve
{Z}(\sigma)$, respectively. Finally, by combining the heteroclinic and
homoclinic parts of the evolution of bifurcations, we obtain the desired
result, thereby completing the proof.
\end{proof}
\end{thm}

As in the case of Theorem 1, the proof of Theorem 2 yields more information
than provided in the statement of the result. For example there are compound
\textquotedblleft blinking effects\textquotedblright\ observed for the
attractor as $\sigma:=C$ increases caused by sequences of tangent
intersections: first the heteroclinic tangencies between $W^{u}(\boldsymbol{p}%
(\sigma))$ and $W^{s}(\boldsymbol{q}(\sigma))$ together with those between
$W^{u}(\boldsymbol{q}(\sigma))$ and $W^{s}(\boldsymbol{p}(\sigma))$, followed
by the homoclinic tangencies between the stable and unstable manifolds of
$\boldsymbol{p}(\sigma)$ together with those of $\boldsymbol{q}(\sigma)$.
Whence, the analysis in \cite{New} together with the proof of Theorem 2 leads
directly to a verification of the following analog of Corollary 1.

\begin{corollary}
Let $f$ be as in Theorem
$2$. Then the bifurcation values $\sigma_{\ast}^{(1)}$ and
$\sigma_{\ast}^{(2)}$ precipitating the creation of a robust stable
chaotic strange attractors for $\sigma_{\ast}^{(1)}<\sigma<\sigma_{t}^{(2)}$
and $\sigma_{\ast}^{(2)}<\sigma<b_{1}$ resulting from
heteroclinic and homoclinic interactions, respectively, are proceeded by
sequences of bifurcation values $\sigma_{t}^{(1)}<\sigma_{t(1)}^{(1)}%
<\sigma_{t(2)}^{(1)}<\cdots<\sigma_{\ast}^{(1)}$ and $\sigma_{t}%
^{(2)}<\sigma_{t(1)}^{(2)}<\sigma_{t(2)}^{(2)}<\cdots<\sigma_{\ast}^{(2)}$
corresponding to successive tangent intersections between the stable and
unstable manifolds $W^{s}\left(  \boldsymbol{p}(\sigma)\right)  $ and
$W^{u}\left(  \boldsymbol{q}(\sigma)\right)  $ and $W^{s}\left(
\boldsymbol{q}(\sigma)\right)  $ and $W^{u}\left(  \boldsymbol{p}
(\sigma)\right)  $, and $W^{s}\left(  \boldsymbol{p}(\sigma)\right)$
and $W^{u}\left(  \boldsymbol{q}(\sigma)\right)  $ and
$W^{s}\left(  \boldsymbol{q}(\sigma)\right)  $ and $W^{u}\left(
\boldsymbol{p}(\sigma)\right)  $ as $\sigma$ increases.
\end{corollary}

It is worth mentioning that the heteroclinic-homoclinic bifurcation is apt to
experience dynamical crises for smaller parameter values than the homoclinic
type because there are two unstable manifolds, rather than just one, capable
of interacting with each stable manifold.

\section{Some Higher-Dimensional Generalizations}

There is an extensive array of possible generalizations of Theorems 1 and 2,
some of which might be quite difficult to realize in terms of discrete dynamical
systems in three or more dimensions. For example, suppose that instead of an
invariant attracting closed curve, we have an attracting invariant torus in
$\mathbb{R}^{3}$ on which the restricted dynamics is ergodic. In addition,
assume there is a single saddle point with a 2-dimensional stable manifold and
the analog of a slice set to which the torus tends to as the parameter of
choice increases. Then one would expect an extremely interesting and complex
analog of the 2-dimensional bifurcation described in Theorem 1. However,
finding a relatively simple smooth map of $\mathbb{R}^{3}$ satisfying these
properties is rather difficult, so we shall confine ourselves to a couple of
much simpler examples.

\subsection{A simple 3-dimensional generalization}

\begin{figure}[htbp]
\centering
\includegraphics[width=0.32\textwidth]{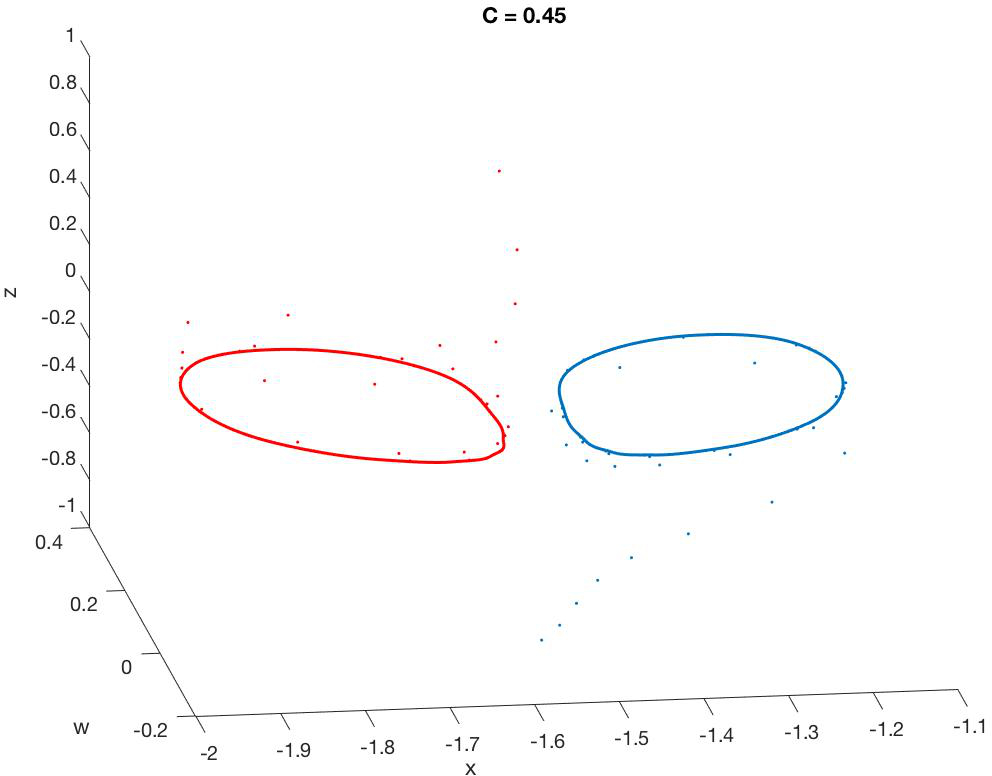}
\includegraphics[width=0.32\textwidth]{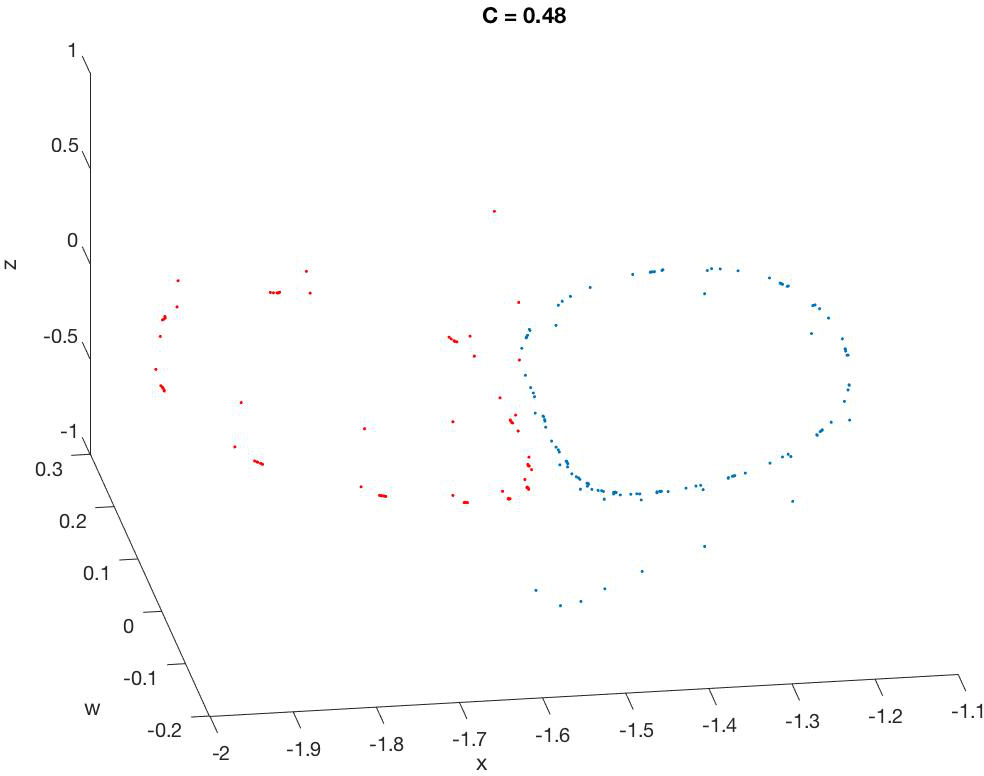}
\includegraphics[width=0.32\textwidth]{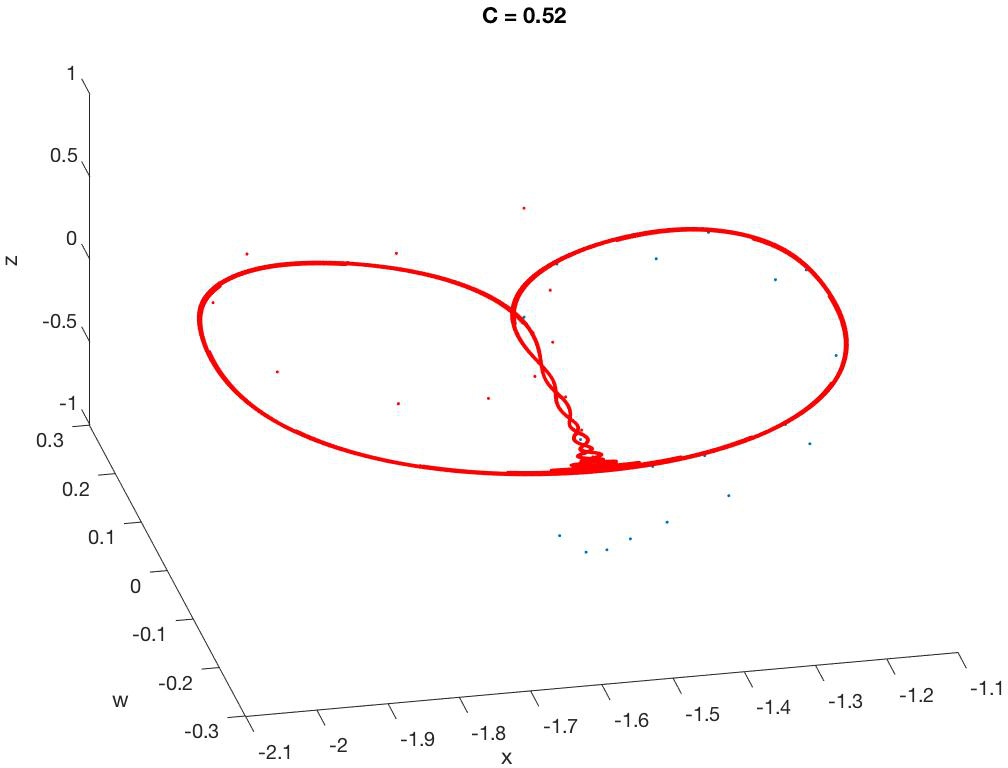}
\caption{Bifurcation evolution for the map $E$ \eqref{e5}.}%
\label{gen1}%
\end{figure}

The first generalization is a more or less trivial extension of Gilet's planar
map; namely%

\begin{equation}
E(x,y,z;\sigma):=\left(  x-\sigma\Psi^{\prime}(x)y,\mu(y+\Psi(x)),0.8z\right)
. \label{e5}%
\end{equation}
To see this, we note that the fixed points of $E$ are $(x_{\ast},y_{\ast},0)$,
where $(x_{\ast},y_{\ast})$ are the fixed points of Gilet's planar map
\eqref{e2} and that the $x$- and $y$-coordinate maps are independent of $z$.
Therefore, the fixed points comprise denumerably many hyperbolic points of the
form $\left(  x_{k},\mu(1-\mu)^{-1}\Psi(x_{k}),0\right)  $, with $\Psi
^{\prime}(x_{k})=0$, each having a 2-dimensional stable and 1- dimensional
unstable manifold, together with a denumerable set of hyperbolic fixed points
$\left(  \tilde{x}_{m},0,0\right)  $, with $\Psi(x_{m})=0$, each of which has
a 3-dimensional stable manifold that bifurcates into a 1-dimensional stable
manifold with a 2-dimensional unstable as the parameter $\sigma$ increases.
These later fixed points, associated with the Neimark--Sacker bifurcations of
Gilet's map, are on invariant lines $x=x_{m},y=0$, which bifurcate into
invariant attracting cylinders of the form $C(\sigma)\times\mathbb{R}$, where
the curves $C(\sigma)$ are as defined in Theorem 1. A simulation of $E$ for
increasing $\sigma$ is shown in Fig. \ref{gen1}.

\subsection{A more complex 3-dimensional generalization}

\begin{figure}[htbp]
\centering
\includegraphics[width=0.49\textwidth]{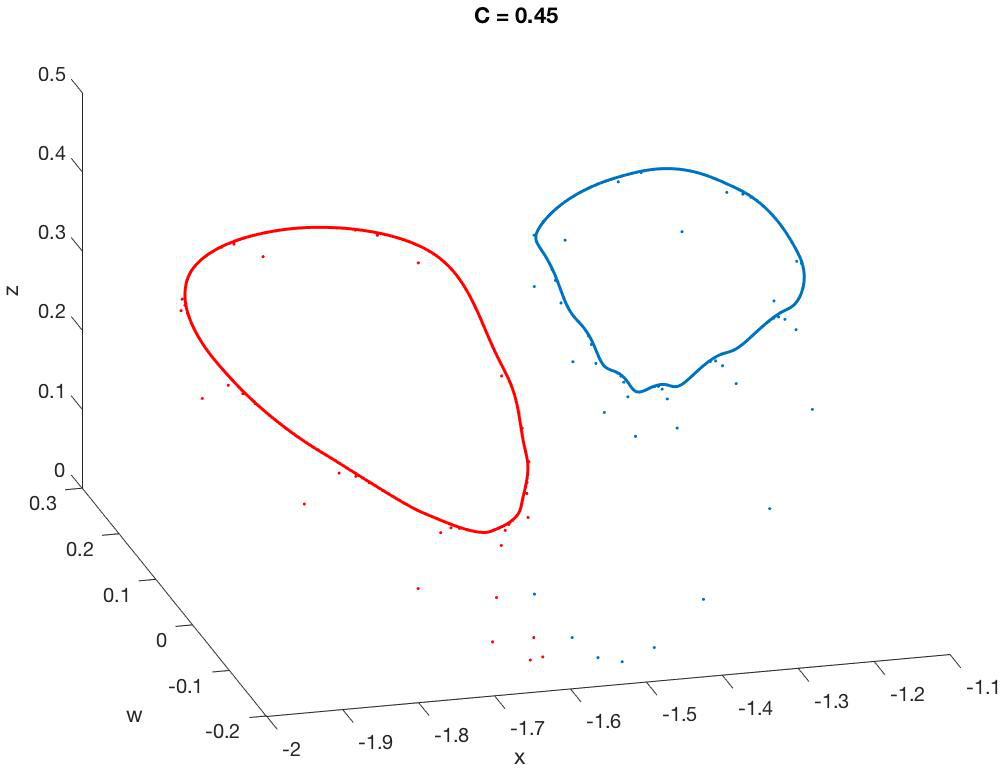}
\includegraphics[width=0.49\textwidth]{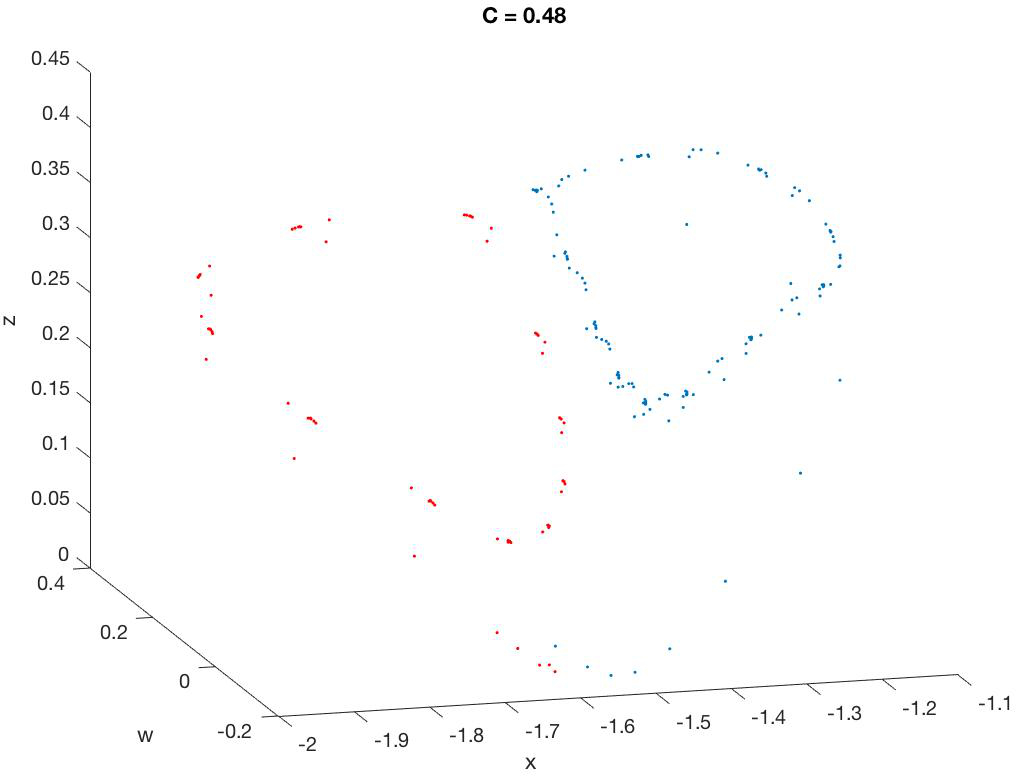}
\includegraphics[width=0.49\textwidth]{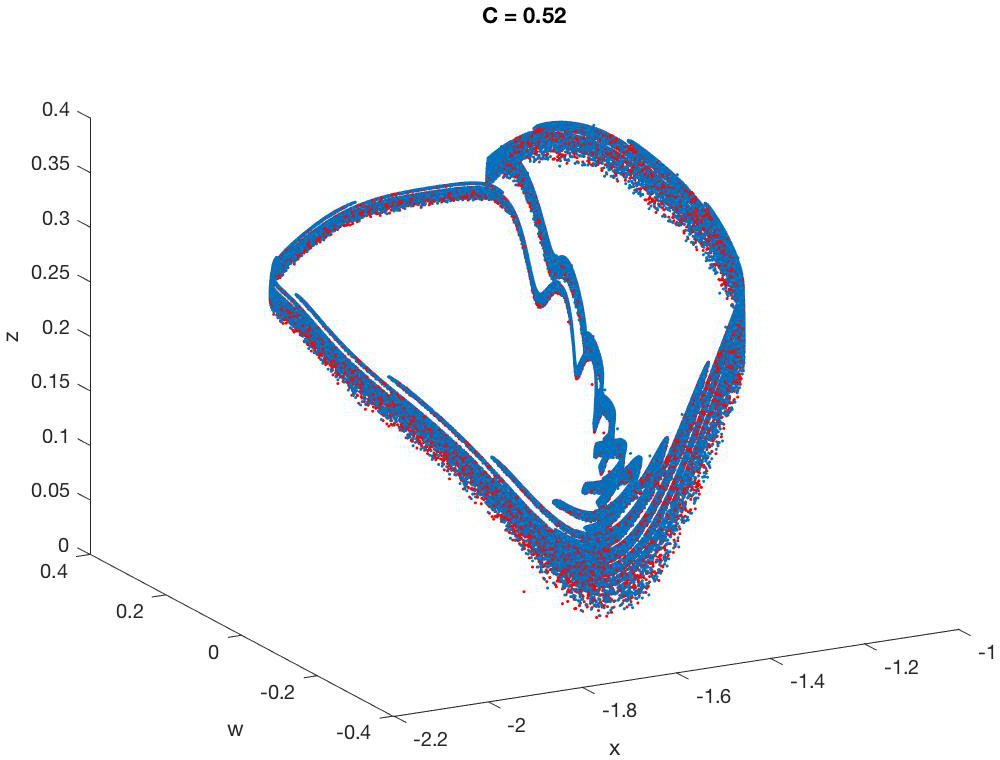}
\includegraphics[width=0.49\textwidth]{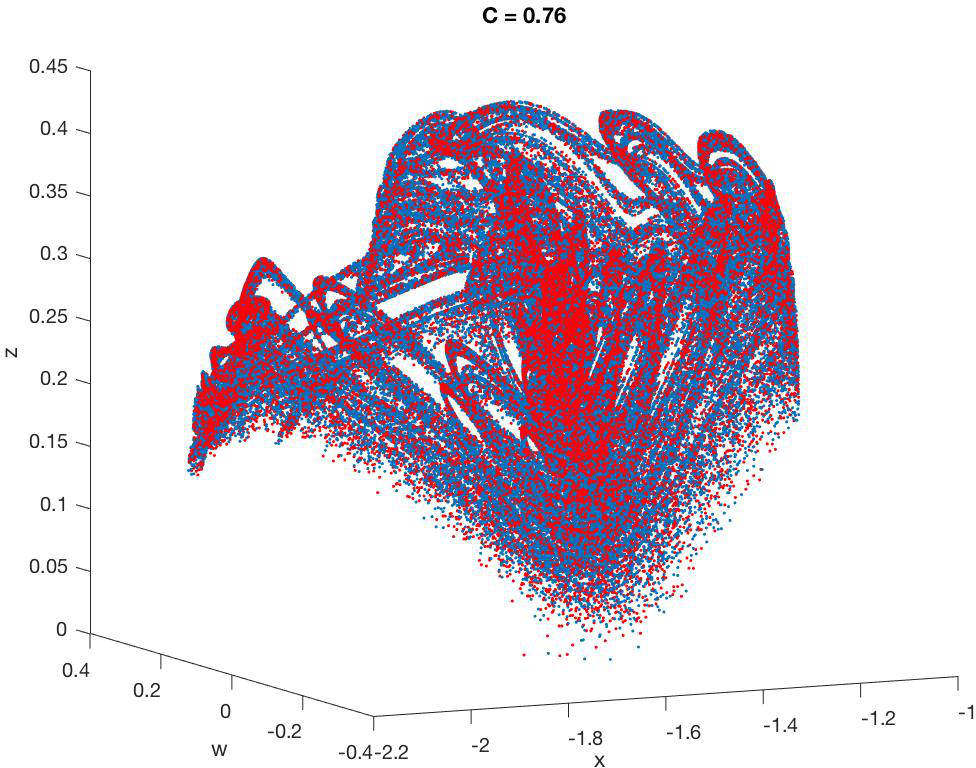}
\caption{Bifurcation evolution for the map $\hat{E}$ \eqref{e6}.}%
\label{gen2}%
\end{figure}

We can create a more interesting extension of Gilet's map by adding a
dependence of the $z$-coordinate map on the other variables; for example, as
in
\begin{equation}
\hat{E}(x,y,z;\sigma):=\left(  x-\sigma\Psi^{\prime}(x)y,\mu(y+\Psi
(x)),0.8z+0.1\sin^{2}(z+\Psi^{\prime}(x)\right)  . \label{e6}%
\end{equation}
It is easy to verify that the fixed points of \eqref{e6} are as follows: There
are, as for \eqref{e5}, denumerably many fixed points of the type $\left(
x_{k},\mu(1-\mu)^{-1}\Psi(x_{k}),0\right)  $, with $\Psi^{\prime}(x_{k})=0$,
each having a 2-dimensional stable and 1-dimensional unstable manifold. There
are also denumerably many fixed points of the form $\left(  \tilde{x}%
_{m},0,z_{m}\right)  $, with $\Psi(x_{m})=0$ and $z_{m}$ the unique solution
of $2z_{m}=\sin^{2}\left(  z_{m}+\Psi^{\prime}(x_{m})\right)  $. These latter
fixed points start as sinks and bifurcate into hyperbolic fixed points each
with a 2-dimensional unstable manifold and a 1-dimensional stable manifold
parallel to the $z$-axis. Moreover, the cylinders $C(\sigma)\times\mathbb{R}$
are attracting and invariant for sufficiently large values of the parameter
$\sigma$. The main difference between this extension and \eqref{e5} is that
the unstable manifolds for the fixed points $\left(  x_{k},\mu(1-\mu)^{-1}%
\Psi(x_{k}),0\right)  $ do not remain in the $x,y$-plane as they wrap around
the cylinder $C(\sigma)\times\mathbb{R}$, which leads to the more complex
bifurcation behavior shown in Fig. \ref{gen2}.

\section{Examples}

Here we shall present examples of the bifurcation behavior described in
Theorem 1 and Theorem 2, mainly through illustrations based on simulations and
a bit of analysis. The examples, which comprise a homoclinic and
heteroclinic-homoclinic type, are to be planar maps based upon Gilet's model
and a minor modification thereof.

\subsection{Symmetric homoclinic bifurcations for Gilet's map}

\begin{figure}[htbp]
\centering
\includegraphics[width = 0.45\textwidth]{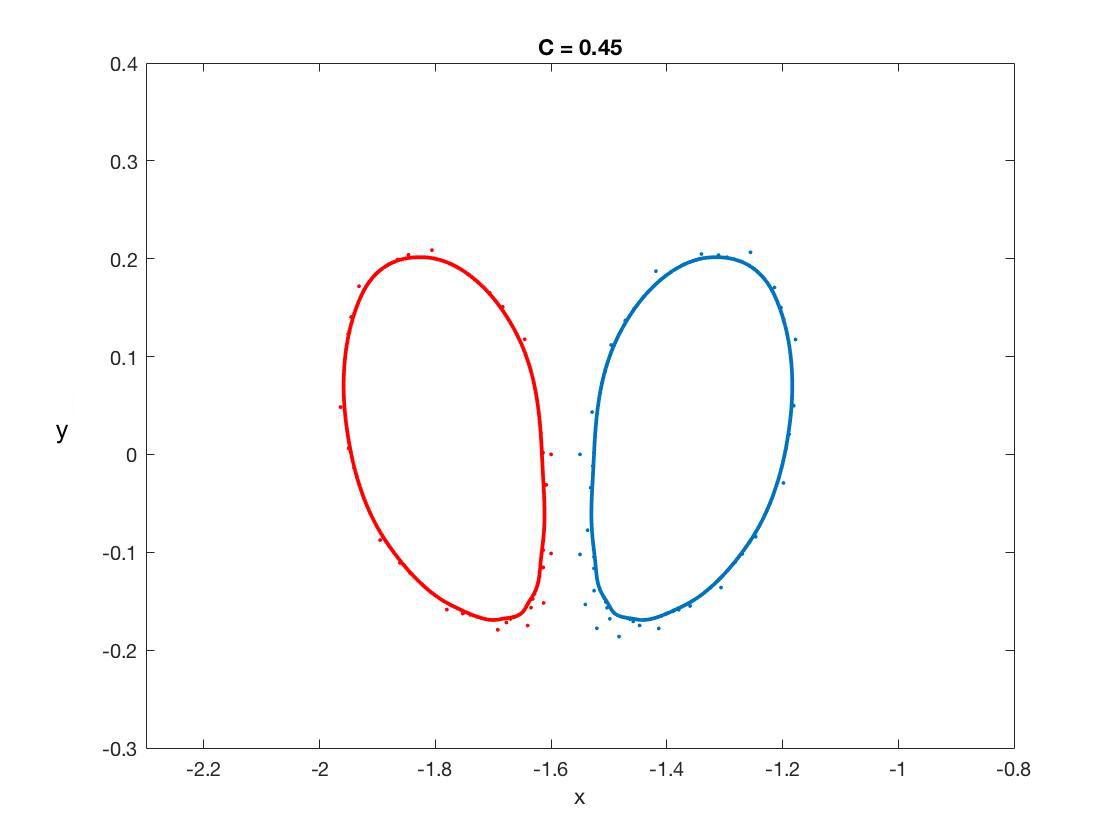}
\includegraphics[width = 0.45\textwidth]{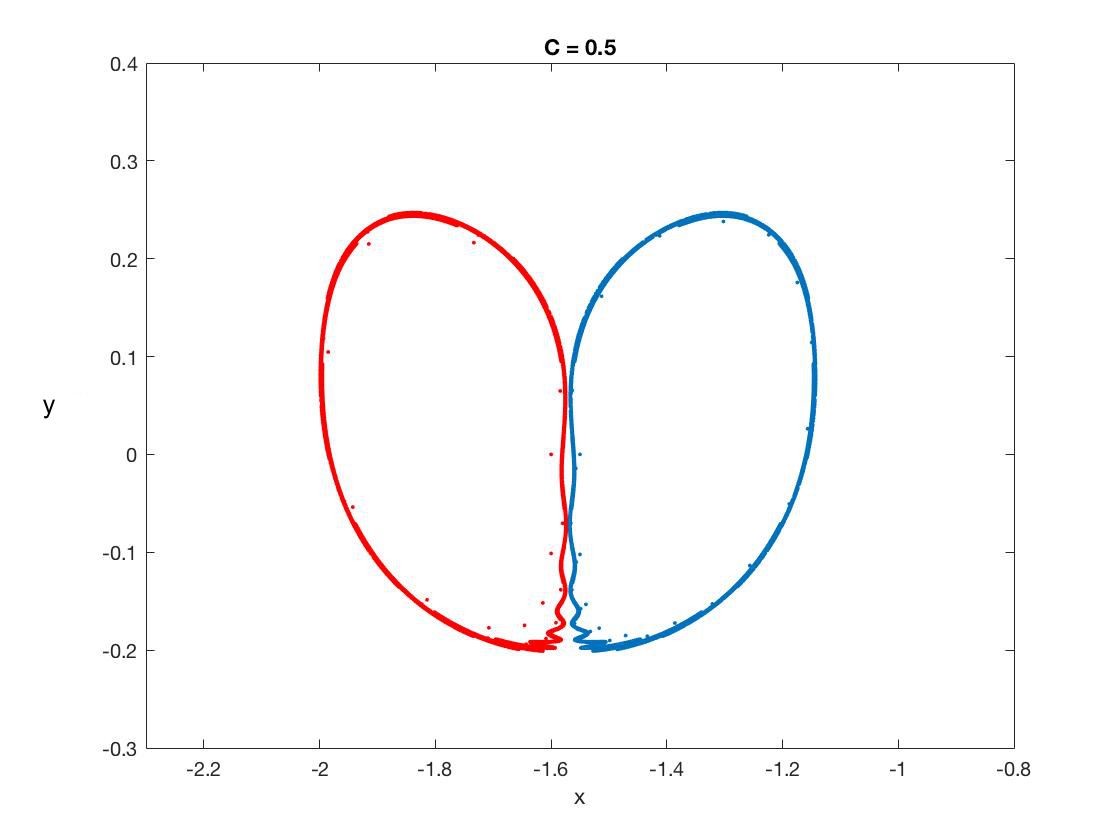}
\includegraphics[width = 0.45\textwidth]{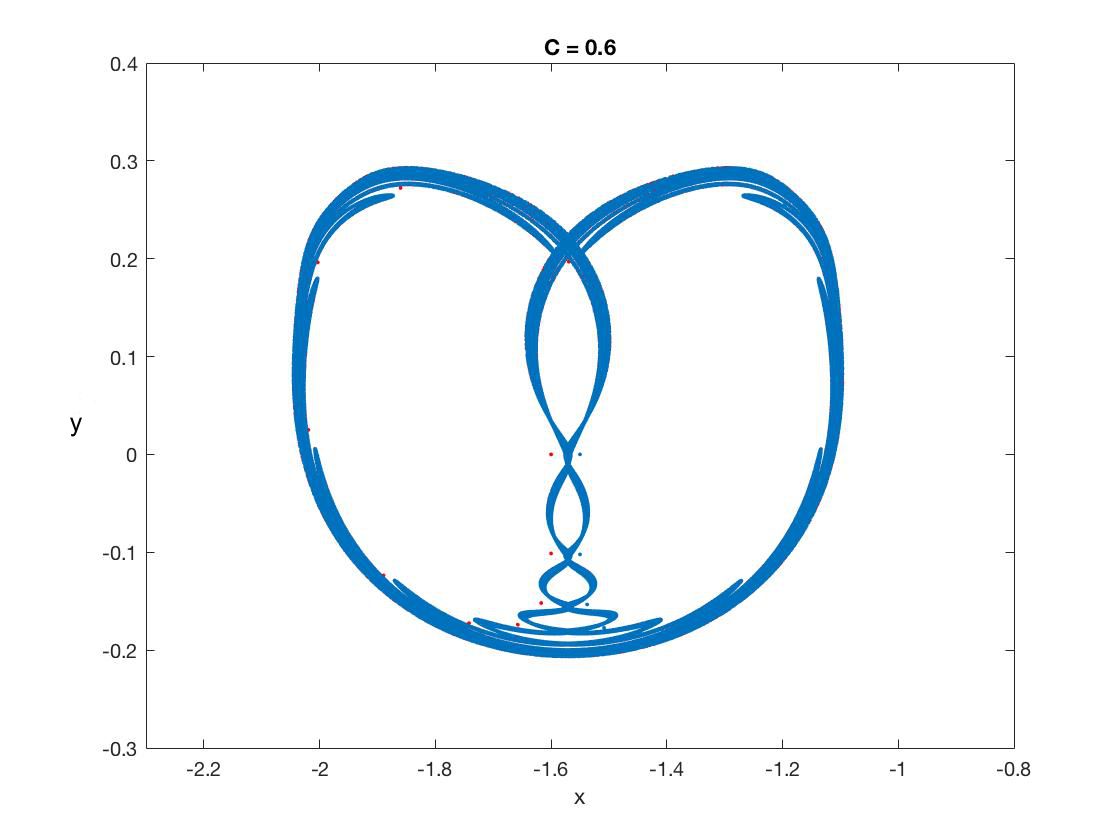}
\includegraphics[width = 0.45\textwidth]{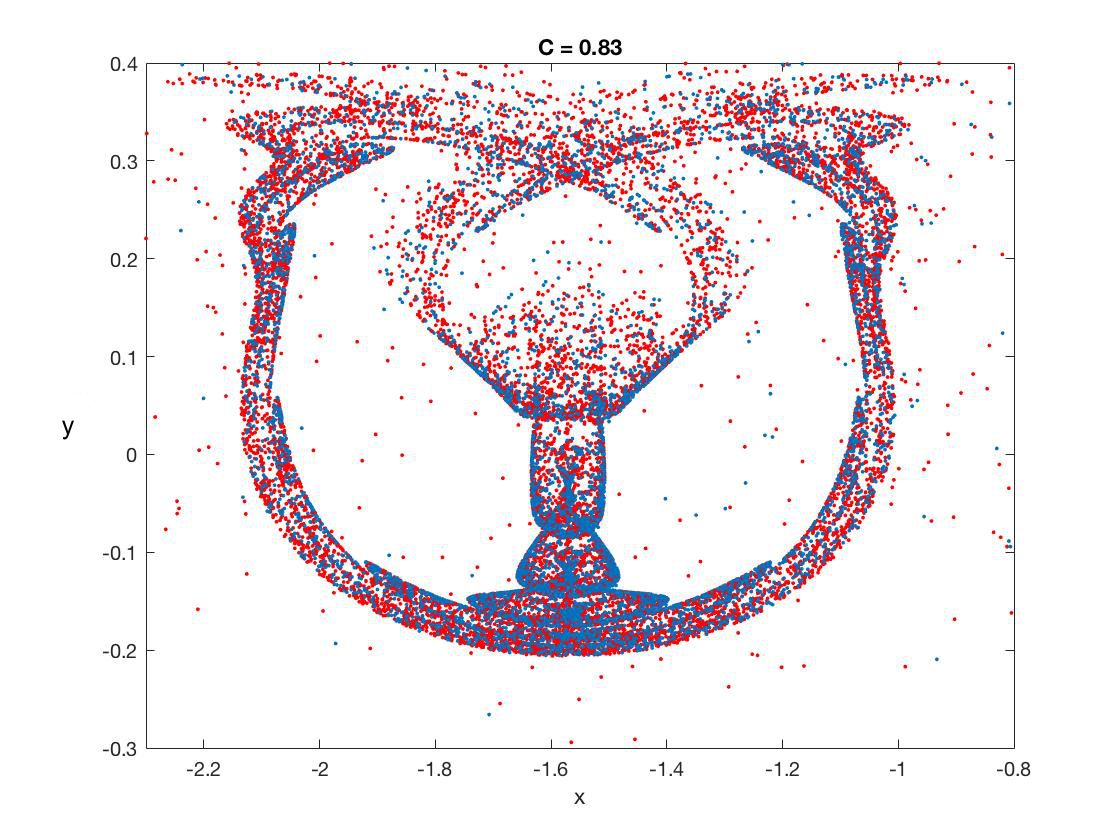}
\caption{Bifurcation evolution for the map $F$ \eqref{e7}.}
\label{Fig: Homoclinic}
\end{figure}

It turns out that homoclinic bifurcations of the type described in Theorem 1
can be very effectively illustrated by considering a pair of contiguous
symmetric cells (containing symmetric invariant closed curves with source
point centers along the $x$-axis) for Gilet's map%

\begin{equation}
F(x,y;\sigma):=\left(  x-\sigma\Psi^{\prime}(x)y,\mu(y+\Psi(x))\right)  ,
\label{e7}%
\end{equation}
where $\mu$ is fixed and $\sigma=C$ is varied. The symmetric cells are chosen
to be the ones symmetric about the stable manifold corresponding to the line
$x\simeq-1.57$, with symmetric center sources (approximately) at $\left(
-1.79,0\right)  $ and $\left(  -1.35,0\right)  $ as shown in Fig. \ref{Fig: Homoclinic}.
It should be noted that, in contrast to the description in Theorem
1, the saddle point on the stable manifold is below rather than above the
symmetric invariant closed attracting curves, which means that the assumptions
remain the same modulo a reflection in the $x$-axis.

Now, it is a straightforward but rather tedious matter to verify that all of
the hypotheses of Theorem 1 hold (modulo the reflection mentioned above), but
all of the assumptions are illustrated quite clearly in Fig. \ref{Fig: Homoclinic}
save the existence of the orientation reversing slice regions. Consequently,
we shall restrict our analysis to the identification of these regions for this
example. A simple calculation shows that the saddle point of interest in this
example is fixed at $\boldsymbol{p}=\left(  \hat{x},\mu(1-\mu)^{-1}\Psi
(\hat{x})\right)  \simeq\left(  -1.57,\mu(1-\mu)^{-1}\Psi(-1.57)\right)
\simeq\left(  -1.57,-0.206\mu(1-\mu)^{-1}\right)  .$ We also compute that
\begin{equation*}
\det F^{\prime}=\mu\left[  1-\sigma\left(  y\Psi^{\prime\prime}(x)-(\Psi
^{\prime}(x))^{2}\right)  \right]  ,
\end{equation*}
so it follows from the fact that $\Psi^{\prime\prime}(x)$ is positive $\left(
\simeq\Psi^{\prime\prime}(\hat{x})\right)  $ in a thin vertical strip centered
at $x=\hat{x}$. Hence, it follows that $\det F^{\prime}$ changes from positive
to negative along the stable manifold as $y$ changes from negative to
positive, which signals a reversal of orientation.

The change in orientation also produces the slice regions $Z$ described in
(A2). To see this, it follows from symmetry that it is enough to describe the
slice to the right of the stable manifold $x=\hat{x}$. Of course, the left
boundary curve for this slice is just a portion of the vertical line
$x=\hat{x}$ for $y\geq\tilde{y}>0.$ The right-hand boundary curve for the
slice is just
\begin{equation*}
y=\frac{x-\hat{x}}{\sigma\Psi^{\prime}(x)}%
\end{equation*}
for $x>\hat{x}$, which has a cusp at $\left(  \hat{x},\tilde{y}\right)  $,
where
\begin{equation*}
\tilde{y}=\lim_{x\downarrow\hat{x}}\left(  \frac{x-\hat{x}}{\sigma\Psi
^{\prime}(x)}\right)  .
\end{equation*}
Now it is straightforward to verify that this slice region has the properties
described in (A2) and (A3) of Theorem1, as does the symmetric slice to the
left of $x=\hat{x}.$

\subsection{Heteroclinic-homoclinic bifurcations for modified Gilet's map}

\begin{figure}[htbp]
\centering
\includegraphics[width=0.32\textwidth]{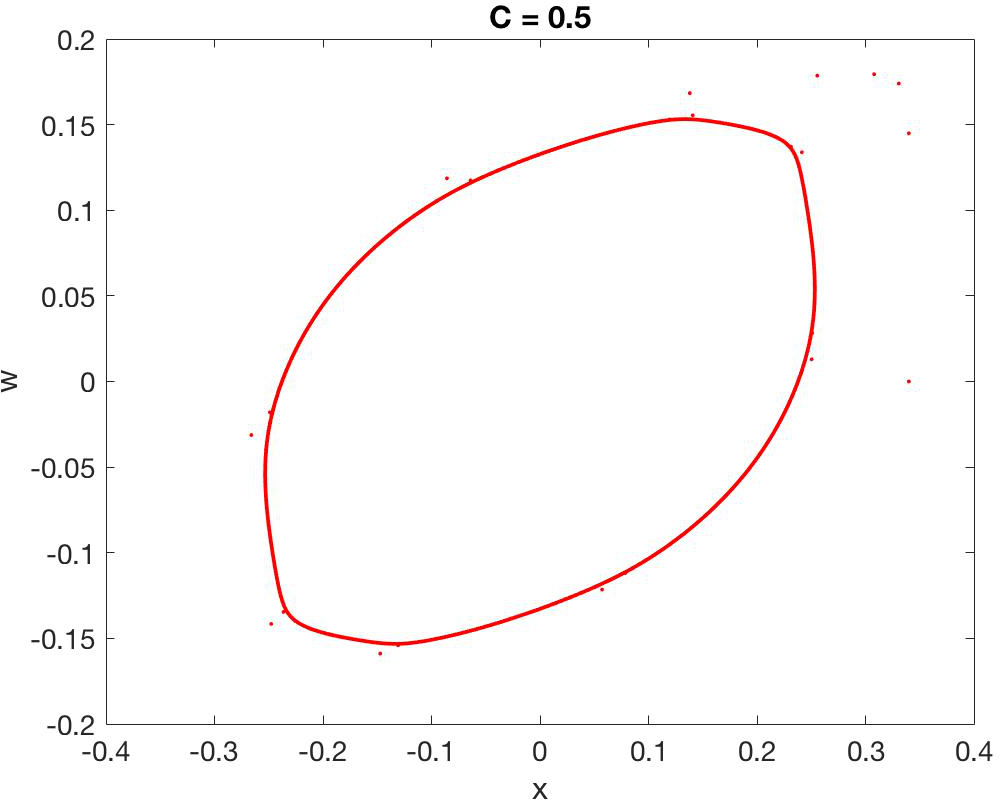}
\includegraphics[width=0.32\textwidth]{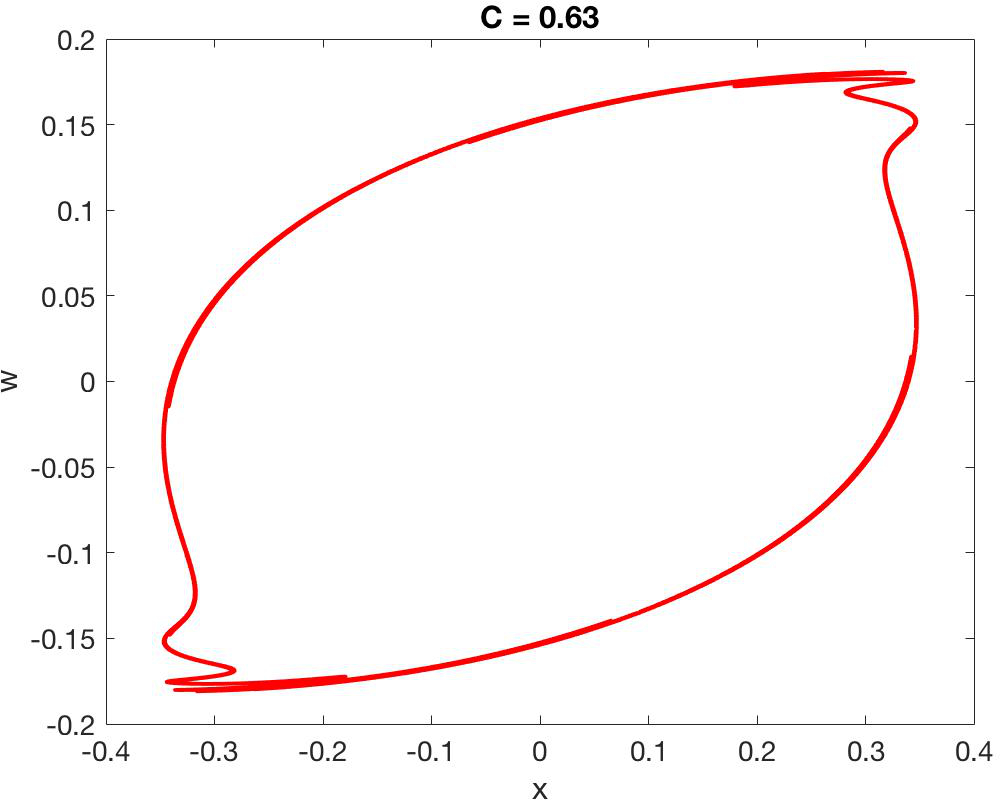}
\includegraphics[width=0.32\textwidth]{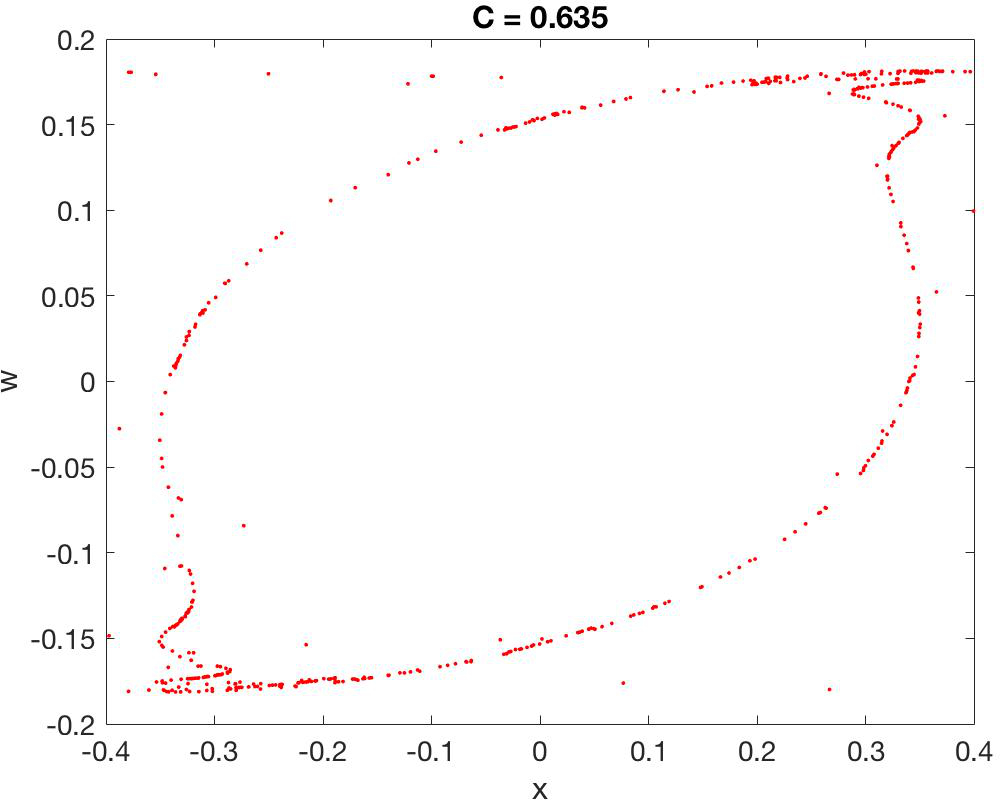}
\caption{Bifurcation evolution for the map $\tilde{F}$ \eqref{e8}.}%
\label{hh1}%
\end{figure}

To illustrate the bifurcations in Theorem 2, we consider the following
modification of Gilet's map: It is the map $\tilde{F}:\mathbb{R}%
^{2}\rightarrow\mathbb{R}^{2}$ defined as follows:%
\begin{equation}
\tilde{F}(x,y;\sigma):=\left\{
\begin{array}
[c]{cc}%
F(x,y;\sigma), & \left\vert x\right\vert \leq0.4\\
\left(  x-[1+5(x+0.4)]\sigma\tilde{\Psi}^{\prime}(x)y,\mu(y+\tilde{\Psi
}(x))\right)  , & -0.6<x<-0.4\\
\left(  x-[1-5(x-0.4)]\sigma\tilde{\Psi}^{\prime}(x)y,\mu(y+\tilde{\Psi
}(x))\right)  , & 0.4<x<0.6\\
\left(  x,\mu(y+\tilde{\Psi}(x))\right)  , & \left\vert x\right\vert \geq0.6
\end{array}
\right.  , \label{e8}%
\end{equation}
where
\begin{equation*}
\tilde{\Psi}(x):=\left\{
\begin{array}
[c]{cc}%
\left(  10x+3\right)  \Psi(0.4), & -0.6<x<-0.4\\
(10x-3)\Psi(0.4), & 0.4<x<0.6\\
-3\Psi(0.4), & x\leq-0.6\\
3\Psi(0.4), & 0.6\leq x
\end{array}
\right.  .
\end{equation*}
We note that the map \eqref{e8} is continuous everywhere and smooth except along
the lines $x=\pm0.4,\pm0.6$. It turns out that the vertical lines along which
it fails to be $C^{\infty}$ do not alter the qualitative nature of the
bifurcation evolution as described in Theorem 2, which is evident from
Fig. \ref{hh1}. Once again, we note that the hypotheses of Theorem 2 can be
readily checked, but all save the existence of the $Z$ slices are clearly
shown in the simulations in Fig. \ref{hh1}. The existence of the orientation
reversing slices $Z(\sigma)$ and $\breve{Z}(\sigma)$ can be verified in a
manner completely analogous to to that used in the previous example.

\section{Concluding Remarks}

In our paper \cite{RB1}, we proved that Gilet's walking droplet model
\cite{Gil} develops Neimark--Sacker bifurcations generating invariant
attracting closed Jordan curves (topological circles) as either one of the parameters is
increased. We also saw that the diameter of these circles increases with
either increasing parameter, which ultimately gives rise to new types of
bifurcations arising from the interactions of stable manifolds with unstable
manifolds of saddle points winding around the expanding circles. The
investigation in this paper comprises an in-depth analysis of two variants -
one purely homoclinic and the other a combination of heteroclinic and
homoclinic interactions of unstable and stable manifolds of saddle points - of
these new bifurcations as the original interaction parameter $C$ (identified
with $\sigma$) is varied while the damping parameter $\mu$ is fixed. In
addition to our analysis of these bifurcations, we showed by examples how
these dynamical phenomena can be extended to higher dimensions.

There are several research directions related to this work that we intend to
pursue in the near future. First among these is a related new bifurcation
somewhat like those studied here, but based on diffeomorphisms that do not
include the orientation reversing, noninjective slice regions in Gilet's
model. Naturally, we also intend to study the bifurcations in models of Gilet
type with $C$ fixed and $\mu$ varying, which, for example, appears to exhibit
more striking dynamical crises behavior than the case studied here. Another
line of research, which has strong connections with the quantum aspects of
pilot waves, is the construction of invariant or approximately invariant
measures for the dynamical models such as Gilet's and is something that we
intend to investigate as part of our continuing investigation of the
mathematical aspects of walking droplet phenomena.

\section*{Acknowledgments}

\noindent Discussions with Anatolij Prykarpatski were very helpful in writing this paper.

\end{document}